\documentclass[12pt, a4paper]{article}
\usepackage[pagewise]{lineno}
\usepackage{graphicx}
\usepackage{amssymb}
\usepackage{latexsym, bm}
\usepackage{multicol}
\usepackage{indentfirst}
\usepackage{amssymb,amsfonts}
\usepackage{amsmath}
\usepackage{cases}
\usepackage[noadjust]{cite}
\usepackage[colorlinks,citecolor=red]{hyperref}
\newtheorem{theorem}{Theorem}
\newtheorem{lemma}{Lemma}
\newtheorem{corollary}{Corollary}

\newtheorem{definition}{Definition}

\newtheorem{remark}{Remark}

\usepackage{amssymb}
\numberwithin{equation}{section}
\numberwithin{theorem}{section}
\numberwithin{remark}{section}
\numberwithin{definition}{section}
\numberwithin{lemma}{section}
\numberwithin{corollary}{section}
\numberwithin{proposition}{section}
\numberwithin{notation}{section}
\title{New Liouville type theorems for the stationary Navier-Stokes equations}
\author{Wenke Tan\footnote{tanwenkeybfq@163.com}\\
{\small Key Laboratory of Computing and Stochastic Mathematics (Ministry of Education),}\\
{\small School of Mathematics and Statistics, Hunan Normal University,}\\
{\small Changsha, Hunan 410081, China}\\
}

\date{}

\begin{document}
\maketitle
{\bf Abstract:} We mainly research the Liouville type problem for the stationary
Navier-Stokes equations (including the fractional case) in $\mathbb{R}^3$. We first establish a new formula for the Dirichlet integral of solutions and show that the globally defined quantity $\int_{\mathbb{R}^3}|\nabla u|^2dx$ is completely determined by the information of the solution $u$ at the origin in frequency space. From this character, we show some new Liouville type theorems for solutions of the stationary Navier-Stokes equations. Then we extend the obtained results for classical stationary Navier-Stokes equations to the stationary fractional Navier-Stokes equations for $\frac{1}{2}\leq s<1$, especially, we solve the Liouville type problem for $s=\frac{5}{6}$.

\medskip

\section{Introduction}
In this paper, we consider the homogeneous stationary Navier-Stokes equations in the whole space $\mathbb{R}^3$:
\begin{equation}\label{NS}
 \left\{\begin{array}{ll}
-\Delta u+u\cdot\nabla u+\nabla P=0\quad in \quad\mathbb{R}^3,\\
\nabla\cdot u=0\quad in\quad\mathbb{R}^3,
\end{array}\right.
\end{equation}
where the unknowns $u(x)=(u_1(x),u_2(x),u_3(x))$ and $P(x)$ denote the velocity vector and the
scalar pressure at the point $x=(x_1,x_2,x_3)\in\mathbb{R}^3$, respectively. Generally, we deal with solutions $u$ of \eqref{NS} in the class of the finite Dirichlet integral
\begin{align}\label{D}
D(u)=\int_{\mathbb{R}^3}|\nabla u|^2dx<+\infty
\end{align} with the homogeneous condition at infinity
\begin{align}\label{far}
\lim_{|x|\to\infty}|u(x)|=0.
\end{align}
It is well known that due to the pioneering work of Leray \cite{Leray}, it has been an open problem whether $u\equiv0$ is the
only solution of \eqref{NS} under conditions \eqref{D} and \eqref{far}. This is a famous Liouville type statement on the stationary Navier-Stokes equations.

For the 2D Navier-Stokes equations, Gilbarg and Weinberger proved the Liouville-type theorem in \cite{G-W}. The authors used the idea that the vorticity $\omega=\nabla\times u$ is a scalar satisfying an elliptic equation that enables one to apply the maximum principle. This, together with another result showing that $\omega\to 0$ at infinity, implies that $\omega=0$. From this, the authors deduced that $u$ and $P$ are constant in $\mathbb{R}^2$. The same approach as in \cite{G-W} fails in the 3D case, mainly due to the more complicated form of the equation for vorticity. Regarding the 3D case, one of the best attempts made to solve the above or related problems was presented by Galdi in \cite{Galdi}, where he showed that if $u\in L^{\frac{9}{2}}(\mathbb{R}^3)$, then it holds that $u\equiv0$. In the case of the $n$-dimensional Navier-Stokes system with $n\geq 4$, the problem was resolved by Galdi in \cite{Galdi} by similar argument to the case of $u\in L^{\frac{9}{2}}(\mathbb{R}^3)$. It is worth pointing out that in the $n$-dimension with $n\geq 5$, the condition \eqref{far} plays a key role in proving Liouville theorems and the Liouville problem without \eqref{far} remains open. Chae and Wolf \cite{C-W} gave a logarithmic improvement of Galdi's
result by assuming that
\begin{align*}
\int_{\mathbb{R}^3}|u|^\frac{9}{2}\{\ln(2+|u|^{-1})\}^{-1}dx<\infty.
\end{align*} Also, Chae \cite{C} showed that the condition
\begin{align*}
\Delta u\in L^\frac{6}{5}(\mathbb{R}^3)
\end{align*} is sufficient for $u=0$ in $\mathbb{R}^3$. He emphasizes that the norm $\Delta u$ in $L^\frac{6}{5}(\mathbb{R}^3)$ corresponds to that of $\nabla u$ in $L^2(\mathbb{R}^3)$ at the level of scaling and that there is no mutual implication relation between their results \cite{C} and \cite{Galdi}. Kozono, Terasawa and Wakasugi proved in \cite{Kozono} that $u=0$ if the vorticity satisfies
\begin{align*}
\lim_{|x|\to \infty}|x|^\frac{5}{3}|\omega(x)|\leq(\delta D(u))^\frac{1}{3}
\end{align*} or the velocity satisfies
\begin{align*}
||u||_{L^{\frac{9}{2},\infty}(\mathbb{R}^3)}\leq (\delta D(u))^\frac{1}{3}
\end{align*} for a small constant $\delta$. Then, the restriction imposed on the norm $||u||_{L^{\frac{9}{2},\infty}(\mathbb{R}^3)}$ in \cite{Kozono} was relaxed by Seregin and Wang in \cite{S-W}. For other values of the parameter
$p\neq\frac{9}{2}$ in the Lebesgue spaces $L^p(\mathbb{R}^3)$.  Chamorro, Jarr\'in and Lemari\'e-Rieusset \cite{C-J-L} prove that $u\in L^p(\mathbb{R}^3)$ with $3\leq p<\frac{9}{2}$
implies $u=0$. When $\frac{9}{2}< p<6$, the authors in \cite{C-J-L} also showed that $u=0$ if $u\in L^p(\mathbb{R}^3)\cap \dot{B}^{\frac{3}{p}-\frac{3}{2}}_{\infty,\infty}$, where $\dot{B}^{\frac{3}{p}-\frac{3}{2}}_{\infty,\infty}$ is a homogeneous Besov space. Moreover, for the value $p=6$, an interesting result of G. Seregin given
in \cite{S} shows that this problem is solved in the space $u\in BMO^{-1}(\mathbb{R}^3)\cap L^6(\mathbb{R}^3)$. Jarr\'in \cite{J} established Liouville type theorem for $u\in L^p(\mathbb{R}^3)$ with $\frac{3}{2}<p<3$ and $u\in L^p(\mathbb{R}^3)\cap \dot{H}^{-1}(\mathbb{R}^3)$ with $\frac{9}{2}<p<\infty$. It is worth pointing out that these results in \cite{C-J-L,S,J} do not suppose the condition \eqref{D} and \eqref{far}. Recently, some results with some local conditions were established by several authors. Seregin \cite{S2} proved that $u=0$ if
\begin{align*}
\sup_{R>0}R^\beta(\frac{1}{|B(R)|}\int_{B(R)}|u|^qdx)^\frac{1}{q}<\infty
\end{align*} for $\frac{3}{2}<q<3$ and $\beta>\frac{6q-3}{8q-6}$. Moreover, he also arrived at the same conclusion under the
assumption
\begin{align*}
\sup_{R>0}R^{\gamma-\frac{3}{s}}||u||_{L^{s,\infty}(B(R))}<\infty,
\end{align*}where $2<s\leq3$ and $\gamma>\frac{4q-3}{6q-6}$.
 Seregin and Wang \cite{S-W} proved the vanishing of $u$ assuming either for $3<q<\infty$, $3\leq l\leq\infty$ (or $q=l=3$),
 \begin{align*}
 \liminf_{R\to\infty}R^{\frac{2}{3}-\frac{2}{q}}||u||_{L^{q,l}(B(R)\setminus B(\frac{R}{2}))}\leq \delta D(u)^\frac{3}{2},
 \end{align*}with $\delta$ small constant, or for $\frac{12}{5}<q<3$, $1\leq l\leq\infty$, $\gamma>\frac{1}{3}+\frac{1}{q}$,
\begin{align*}
 \liminf_{R\to\infty}R^{\gamma-\frac{2}{q}}||u||_{L^{q,l}(B(R)\setminus B(\frac{R}{2}))}<\infty.
 \end{align*} Tsai \cite{Tsai} established the Liouville type theorem, using the assumption
\begin{align*}
 \liminf_{R\to\infty}\frac{1}{R}||u||^{3-\delta}_{L^q(R\leq|x|\leq\lambda R)}=0,
 \end{align*} where $q=q(\delta)=\frac{3-\delta}{1-\frac{\delta}{6}}$, for some constants $0\leq\delta\leq1$ and $\lambda>1$, or the assumption
\begin{align*}
 \liminf_{R\to\infty}R^\beta||u||^{3-\delta}_{L^q(R\leq|x|\leq R+R^{1-\alpha})}=0,
 \end{align*}where $q=q(\delta)=\frac{3-\delta}{1-\frac{\delta}{6}}$ and  $\beta=\max\{\frac{(3-\alpha)q^{-1}-2+3\alpha}{2-\delta}, \frac{-1+2\alpha}{3-\delta}\}$, for some constants $0\leq\delta\leq1$ and $\alpha\geq0$. Cho, Neustupa and Yang \cite{C-N-Y} proved that if $u$ satisfies
 \begin{align*}
 \liminf_{R\to\infty}R^{\frac{1}{3}-\frac{2}{q}}||u||_{L^q(B(2R)\setminus B(R))}<\infty\quad for \quad some\quad \frac{3}{2}<q<3
 \end{align*} or
 \begin{align*}
 \liminf_{R\to\infty}R^{\frac{-1}{3}}||u||_{L^3(B(2R)\setminus B(R))}=0.
 \end{align*}Then $u=0$. It is worth pointing out that the results in \cite{C-N-Y} do not need the condition \eqref{D}.

We now turn our attention to the stationary fractional Navier-Stokes equations in $\mathbb{R}^3$.
\begin{equation}\label{FNS}
 \left\{\begin{array}{ll}
(-\Delta)^s u+u\cdot\nabla u+\nabla P=0\quad for \quad s\in(0,1),\\
\nabla\cdot u=0.
\end{array}\right.
\end{equation}For the non-local operator $(-\Delta)^s$, it is defined for any $f\in\mathcal{S}'(\mathbb{R}^3)$ through the
Fourier transform: $(-\Delta)^sf(x)=\mathcal{F}^{-1}(|\xi|^{2s}\mathcal{F}f(\xi))$.

Compared with the classical stationary Navier-Stokes equations, there are very few results on the steady fractional Navier-Stokes equations \eqref{FNS}.
Tang and Yu \cite{TY} studied partial H\"{o}lder regularity of the steady fractional Navier-Stokes equations with $\frac{1}{2}<s<1$, and proved that if $\frac{1}{2}<s<\frac{5}{6}$,
 any steady suitable weak solution is regular away from a relatively closed set with zero $5-6s$-Hausdorff measure and it is regular when $\frac{5}{6}\leq s<1$.
Then, a natural question is: Are there non-identically vanishing smooth solutions $u$ to \eqref{FNS} with $0<s<1$
satisfying the conditions \eqref{far} and
\begin{align}\label{FD}
\int_{\mathbb{R}^3}|(-\Delta)^\frac{s}{2}u|^2dx<\infty.
\end{align} By using the Caffarelli-Sivestre extension to convert the non-local operator $(-\Delta)^s$
on $\mathbb{R}^3$ to the local operator $\bar{\Delta}$ on $\mathbb{R}^4_+$, Yang \cite{Yang} proved that if $\frac{5}{6}\leq s<1$, $u\in L^\frac{9}{2}(\mathbb{R}^3)$ implies $u\equiv0$. Furthermore, by embedding theorem $\dot{H}^\frac{5}{6}(\mathbb{R}^3)\hookrightarrow L^\frac{9}{2}(\mathbb{R}^3)$, Yang \cite{Yang}
showed $u$ is trivial when $s=\frac{5}{6}$. When $0<s<\frac{5}{6}$, this problem is open.

To the author's knowledge, those authors of the previous works in the literature are analyzing the Liouville problem in physical space. There are no related results that consider the
Liouville problem in frequency space. In this paper, we will exploit the localization techniques in frequency space to research the Liouville problem for the solution to the stationary Navier-Stokes equations \eqref{NS}(the stationary fractional Navier-Stokes equations \eqref{FNS}) with assumptions \eqref{D}(\eqref{FD}).

Our analysis starts from the scaling symmetry $(u,P)\rightarrow(u_\lambda(x),P_\lambda(x))$ of stationary Navier-Stokes equations \eqref{NS} for $\lambda>0$, where
 \begin{align*}
 u_\lambda(x)\doteq \lambda u(\lambda x)\quad P_\lambda(x)\doteq \lambda^2 P(\lambda x).
 \end{align*}
 Among other things, this means the quantity $\int_{\mathbb{R}^3}|\nabla u|^2dx$ is subcritical. This character inspires us that the main obstacle to solving the Liouville problem may come from the low-frequency parts of solutions. To exploit the subcritical character of \eqref{D}, we use the localization techniques in frequency space to transform the Dirichlet integral to a new integral which is completely determined by the behavior of the solution $u$ at the origin in frequency space. The first step in doing so is to observe that the high-frequency part of solutions $u^k$ can be taken as a test function to take $L^2$ inner product to equations \eqref{NS}. Then by using the fact $\nabla\cdot u=0$, we can eliminate the mainly high-frequency parts $\int_{\mathbb{R}^3}u_j\partial_ju^k_iu^k_idx$. Secondly, we apply Bony decomposition to localized the "low-low-high" term $\int_{\mathbb{R}^3}\dot{S}_k u_j\partial_j\dot{S}_ku_iu^k_idx$ into the region near the origin in frequency space. The remainder is the "high-low-high" term $\int_{\mathbb{R}^3}u^k_j\partial_j\dot{S}_ku_iu^k_idx$ which is a globally defined quantity. To overcome the main difficulty, we first apply Bony composition to separate the essential high frequency part $\int_{\mathbb{R}^3}\partial_j\dot{S}_ku_i\sum_{l\geq k-1}\dot{\Delta}_lu^k_i\tilde{\dot{\Delta}}u^k_jdx$ from the original integral $\int_{\mathbb{R}^3}u^k_j\partial_j\dot{S}_ku_iu^k_idx$. The key point is that we can introduce a parameter $\theta$ to split this term into two parts
 \begin{align*}
 &\int_{\mathbb{R}^3}\partial_j\dot{S}_ku_i\sum_{k-1\leq l\leq [\theta k]}\dot{\Delta}_lu^k_i\tilde{\dot{\Delta}}_{l}u^k_jdx+\int_{\mathbb{R}^3}\partial_j\dot{S}_ku_i\sum_{l\geq [\theta k]+1}\dot{\Delta}_lu^k_i\tilde{\dot{\Delta}}_{l}u^k_jdx\end{align*} then by choosing a suitable value for $\theta$ (say $\theta=\frac{1}{2}$), we find that the high-frequency part $\int_{\mathbb{R}^3}\partial_j\dot{S}_ku_i\sum_{l\geq [\theta k]+1}\dot{\Delta}_lu^k_i\tilde{\dot{\Delta}}_{l}u^k_jdx$
will vanish as $k\to-\infty$. We thus localize the Dirichlet integral into the region near the origin in frequency space.

For the fractional cases, we apply the same ideas to localize the fractional Dirichlet integral \eqref{FD}. The main difference is that when $\frac{5}{6}<s<1$, \eqref{FD} is also subcritical and we localize it into the low frequency region; when $\frac{1}{2}\leq s<\frac{5}{6}$, \eqref{FD} is supercritical and we localize it into the high frequency region; when $s=\frac{5}{6}$, \eqref{FD} is critical and we solve the Liouville problem. It is worth pointing out that our methods avoiding the use of the Caffarelli-Sivestre extension provide a unified research framework for classical Navier-Stokes equations and fractional Navier-Stokes equations.

Our first result which establishes a new formula for Dirichlet integral of solution to \eqref{NS} now reads:
\begin{theorem}\label{main}
Let $u$ be a smooth solution of \eqref{NS} in the class \eqref{D}. Then the following identity holds
\begin{align}\label{snc}
&\int_{\mathbb{R}^3}|\nabla u|^2dx\\=&-\liminf_{k\to-\infty}\{\int_{\mathbb{R}^3}
\sum_{l=k-1}^{k+2}\sum_{l'=l-2}^{k-1}\dot{S}_{l-2}\dot{S}_{k}u\cdot\nabla\dot{\Delta}_{l'}u\cdot\dot{\Delta}_lu^kdx+\int_{\mathbb{R}^3}\sum_{l=k-3}^{ k}\dot{\Delta}_l\dot{S}_ku\cdot\nabla\dot{S}_ku\cdot\tilde{\dot{\Delta}}_lu^kdx\notag\\&+\int_{\mathbb{R}^3}\sum_{l=k-1}^{[\frac{k}{2}]}\dot{\Delta}_l u^k\cdot\nabla \dot{S}_ku\cdot\tilde{\dot{\Delta}}_lu^kdx\},\notag
\end{align}here we have used the notation that $\tilde{\dot{\Delta}}_lu=\sum_{|l'-l|\leq2}\dot{\Delta}_{l'}u$.
\end{theorem}
\begin{remark}
Notice that
\begin{align*}
&\liminf_{k\to-\infty}supp\mathcal{F}(\sum_{l=k-1}^{k+2}\sum_{l'=l-2}^{k-1}\dot{S}_{l-2}\dot{S}_{k}u\cdot\nabla\dot{\Delta}_{l'}u\cdot\dot{\Delta}_lu^k)=\{0\},\\
&\liminf_{k\to-\infty}supp\mathcal{F}(\sum_{l=k-3}^{k}\dot{\Delta}_l\dot{S}_ku\cdot\nabla\dot{S}_ku\cdot\tilde{\dot{\Delta}}_lu^k)=\{0\},\\
&\liminf_{k\to-\infty}supp\mathcal{F}(\sum_{l=k-1}^{[\frac{k}{2}]}\dot{\Delta}_l u^k\cdot\nabla \dot{S}_ku\cdot\tilde{\dot{\Delta}}_lu^k)=\{0\}.
\end{align*}The formula \eqref{snc} means that the globally defined quantity $\int_{\mathbb{R}^3}|\nabla u|^2dx$ is completely characterized by the local information of $u$ at the origin in frequency space.
\end{remark}
As an application of the above theorem, we have the following uniqueness result on \eqref{NS}, specifically, the condition \eqref{far} is removed.
\begin{corollary}(Liouville type theorem)\label{coro}
Let $u$ be a smooth solution of \eqref{NS} in the class \eqref{D}. Then
$u\equiv0$ if the following condition is valid
\begin{align}
\liminf_{k\to-\infty}2^{-k}||\dot{S}_ku||_{L^\infty(\mathbb{R}^3)}<\infty.\label{sc1}
\end{align}Specifically, we deduce that if $u$ satisfies
\begin{align}
&\liminf_{k\to-\infty}||\dot{S}_ku||_{\dot{B}^{-1}_{\infty,\infty}}<\infty,\label{sc2}\\
\quad or\quad&\liminf_{k\to-\infty}||\hat{u}||_{L^r(B_{2^k})}<\infty \quad for \quad r\geq\frac{3}{2},\label{sc3}
\end{align}then $u\equiv0.$
\end{corollary}

Next, we turn our attention to the fractional cases. We first present the
definition of a weak solution of the equations \eqref{FNS}:
\begin{definition}
We say that the function $u\in \dot{H}^s(\mathbb{R}^3)$ is a weak solution to equations \eqref{FNS} if $u$ is weakly divergence free and the following identity holds
\begin{align*}
\int_{\mathbb{R}^3}(-\Delta)^\frac{s}{2}u\cdot(-\Delta)^\frac{s}{2}\varphi dx=\int_{\mathbb{R}^3}(u\otimes u):\nabla\varphi dx+\int_{\mathbb{R}^3}P\nabla\cdot\varphi dx.
\end{align*} for some $P\in\mathcal{S}'(\mathbb{R}^3)$ and all $\varphi(x)\in (C^\infty_0(\mathbb{R}^3))^3$.
\end{definition}
Our first result for stationary fractional Navier-Stokes equations \eqref{FNS} now reads as: 
\begin{theorem}\label{main1}
Assume $\frac{5}{6}\leq s<1$ and $u$ be a weak solution of \eqref{FNS} in the class \eqref{FD}. Then the following identity holds
\begin{align}\label{snc1}
&\int_{\mathbb{R}^3}|(-\Delta)^\frac{s}{2} u|^2dx\\=&-\liminf_{k\to-\infty}\{\int_{\mathbb{R}^3}
\sum_{l=k-1}^{k+2}\sum_{l'=l-2}^{k-1}\dot{S}_{l-2}\dot{S}_{k}u\cdot\nabla\dot{\Delta}_{l'}u\cdot\dot{\Delta}_lu^kdx+\int_{\mathbb{R}^3}\sum_{l=k-3}^{k}\dot{\Delta}_l\dot{S}_ku\cdot\nabla\dot{S}_ku\cdot\tilde{\dot{\Delta}}_lu^kdx\notag\\&+\int_{\mathbb{R}^3}\sum_{l=k-1}^{[\frac{k}{2}]}\dot{\Delta}_l u^k\cdot\nabla \dot{S}_ku\cdot\tilde{\dot{\Delta}}_lu^kdx\}.\notag
\end{align}
\end{theorem}
\begin{remark}
The formula \eqref{snc1} means that although the nonlocal operator $(-\Delta)^s$ with $\frac{5}{6}\leq s<1$ is considered, we can still localize the globally defined quantity \eqref{FD} into the region near the origin in frequency space. Furthermore, the localization is independent of $s$.
\end{remark}
Based on the above formula, we also obtain the following Liouville type theorem:
\begin{corollary}(Liouville type theorem)\label{coro1}
Assume $\frac{5}{6}\leq s<1$ and $u$ be a weak solution of \eqref{FNS} in the class \eqref{FD}. Then
$u\equiv0$ if the following condition is valid
\begin{align}
\liminf_{k\to-\infty}2^{k(1-2s)}||\dot{S}_ku||_{L^\infty(\mathbb{R}^3)}<\infty.\label{fsc1}
\end{align}Specifically, we deduce that if $u$ satisfies
\begin{align}
&\liminf_{k\to-\infty}||\dot{S}_ku||_{\dot{B}^{1-2s}_{\infty,\infty}}<\infty,\label{fsc2}\\
\quad or\quad&\liminf_{k\to-\infty}||\hat{u}||_{L^r(B_{2^k})}<\infty \quad for \quad r\geq\frac{3}{4-2s},\label{fsc3}
\end{align}then $u\equiv0.$
\end{corollary}
Specifically, when $s=\frac{5}{6}$, we can get the following uniqueness result.
\begin{corollary}(Uniqueness for $s=\frac{5}{6}$)\label{coro2}
Assume $s=\frac{5}{6}$ and $u$ be a weak solution of \eqref{FNS} in the class \eqref{FD}. Then
$u\equiv0$.
\end{corollary}

When $\frac{1}{2}<s<\frac{5}{6}$, in view of the supercritical character of \eqref{FD}, we localize the fractional Dirichlet integral \eqref{FD} into the region with high frequency.
\begin{theorem}\label{main2}
Assume $\frac{1}{2}< s<\frac{5}{6}$ and $u$ be a weak solution of \eqref{FNS} in the class \eqref{FD}. Then it holds that
\begin{align}\label{snc2}
&\int_{\mathbb{R}^3}|(-\Delta)^\frac{s}{2}u|^2dx\\=&\liminf_{k\to\infty}\{\int_{\mathbb{R}^3}
\sum_{l=k-1}^{k+2}\sum_{l'=l-2}^{k-1}\sum_{l''=[\theta k]}^{k-1}\dot{S}_{l-2}\dot{\Delta}_{l''}u\cdot\nabla\dot{\Delta}_{l'}u\cdot\dot{\Delta}_lu^kdx\notag\\
&+\int_{\mathbb{R}^3}\sum_{l=k-3}^{k}\sum_{l''=[\frac{k}{2}]}^{k-1}\dot{\Delta}_l\dot{S}_ku\cdot\nabla\dot{\Delta}_{l''}u\cdot\tilde{\dot{\Delta}}_lu^kdx
+\int_{\mathbb{R}^3}\sum_{l''=[\frac{k}{2}]}^{k-1}\sum_{k-1\leq l}\tilde{\dot{\Delta}}_{l}u^k\cdot\nabla\dot{\Delta}_{l''}u\cdot\dot{\Delta}_lu^kdx\}\notag
\end{align}here $\theta=\frac{4s-2}{3-2s}\in(0,1)$.
\end{theorem}
\begin{corollary}(Liouville type theorem)\label{coro3}
Assume $\frac{1}{2}< s<\frac{5}{6}$ and $u$ be a weak solution of \eqref{FNS} in the class \eqref{FD}. If $u$ satisfies the following condition
\begin{align}\label{+1}
\liminf_{k\to\infty}2^{k(1-2s)}(||\sum_{[\theta k]\leq l\leq k-1}\dot{\Delta}_lu_j||_{L^\infty}+||\sum_{[\frac{k}{2}]\leq l\leq k-1}\dot{\Delta}_lu_i||_{L^\infty})<\infty,
\end{align} then it holds that $u\equiv0$. Specifically, if $u$ satisfies
\begin{align}\label{+3}
\liminf_{k\to\infty}||u^k||_{\dot{B}^{1-2s}_{\infty,\infty}}<\infty,
\end{align}here $u^k=u-\dot{S}_{k}u$. Then it follows $u\equiv0$.
\end{corollary}
\begin{remark}
From the formula $\theta=\frac{4s-2}{3-2s}$, it is clear that $0<\theta\leq\frac{1}{2}$ when $\frac{1}{2}<s\leq\frac{7}{10}$ and $\frac{1}{2}<\theta<1$ when $\frac{7}{10}<s<\frac{5}{6}$. This means that the condition \eqref{+1} can be simplified as $\liminf_{k\to\infty}2^{k(1-2s)}||\sum_{[\theta k]\leq l\leq(k-1)}\dot{\Delta}_lu_j||_{L^\infty}<\infty$ with $\frac{1}{2}<s\leq\frac{7}{10}$ and $\liminf_{k\to\infty}2^{k(1-2s)}||\sum_{[\frac{k}{2}]\leq l\leq(k-1)}\dot{\Delta}_lu_j||_{L^\infty}<\infty$ with $\frac{7}{10}<s<\frac{5}{6}$.
\end{remark}

When $s=\frac{1}{2}$, we can not establish a localized formula for the fractional Dirichlet integral \eqref{FD}. Our result reads as follows:
\begin{theorem}\label{main3}
Assume $s=\frac{1}{2}$ and $u$ be a weak solution of \eqref{FNS} in the class \eqref{FD}. Then it holds that
\begin{align}\label{snc3}
&\int_{\mathbb{R}^3}|(-\Delta)^\frac{1}{4}u|^2dx\\=&\liminf_{k\to\infty}\{\int_{\mathbb{R}^3}
\sum_{l=k-1}^{k+2}\sum_{l'=l-2}^{k-1}\sum_{l''=0}^{k-1}\dot{S}_{l-2}\dot{\Delta}_{l''}u\cdot\nabla\dot{\Delta}_{l'}u\cdot\dot{\Delta}_lu^kdx\notag\\
&+\int_{\mathbb{R}^3}\sum_{l=k-3}^{k}\sum_{l''=[\frac{k}{2}]}^{k-1}\dot{\Delta}_l\dot{S}_ku\cdot\nabla\dot{\Delta}_{l''}u\cdot\tilde{\dot{\Delta}}_lu^kdx
+\int_{\mathbb{R}^3}\sum_{l''=[\frac{k}{2}]}^{k-1}\sum_{k-1\leq l}\tilde{\dot{\Delta}}_{l}u^k\cdot\nabla\dot{\Delta}_{l''}u\cdot\dot{\Delta}_lu^kdx\}.\notag
\end{align}Furthermore, assume that $u$ satisfies the following condition
\begin{align}\label{+2}
||u^0||_{L^\infty}<\infty
\end{align}here $u^0=u-\dot{S}_0u$.
Then it holds that $u\equiv0$.
\end{theorem}
Our paper is organized as follows: In Section 2, we list some notations and recall some
Lemmas about the Littlewood-Paley theory which will be used in the sequel. In Section 3 we present the proofs of results.

\section{Preliminaries}

{\bf Notation:} In this paper, we denote $B_r=\{x\in\mathbb{R}^3:|x|\leq r\}$.  The matrix $\nabla u$ denotes the gradient of $u$ with respect to the $x$ variable, whose
$(i,j)$-th component is given by $(\nabla u)_{i,j}=\partial_ju_i$ with $1\leq i,j\leq3$. Throughout this
paper, $C$ stands for some positive constant, which may vary from line to line. Given a Banach space $X$, we denote its norm by $||\cdot||_{X}$.
We use $\mathcal{S}(\mathbb{R})^3$ and $\mathcal{S}'(\mathbb{R}^3)$ to denote Schwartz functions and the tempered distributions spaces on $\mathbb{R}^3$, respectively.

Next, we will recall some well-known facts about the Littlewood-Paley decomposition. Firstly, let us recall that for every $f\in\mathcal{S}'(\mathbb{R}^3)$
the Fourier transform of $f$ is defined by
\begin{align*}
(\mathcal{F}f)(\xi)=\hat{f}(\xi)=(2\pi)^{-\frac{3}{2}}\int_{\mathbb{R}^3}e^{-ix\cdot\xi}f(x)dx\quad for\quad\xi\in\mathbb{R}^3
\end{align*}The inverse Fourier transform of $g\in\mathcal{S}'(\mathbb{R}^3)$ is given by
\begin{align*}
(\mathcal{F}^{-1}g)(x)=\breve{g}(x)=(2\pi)^{-\frac{3}{2}}\int_{\mathbb{R}^3}e^{ix\cdot\xi}g(\xi)dx\quad for \quad x\in\mathbb{R}^3.
\end{align*}By using Fourier transform, we can define homogeneous Sobolev norm $||\cdot||_{\dot{H}^s}$ with $s\in\mathbb{R}$ as
\begin{align*}
||f||_{\dot{H}^s(\mathbb{R}^3)}=(\int_{\mathbb{R}^3}|\xi|^{2s}|\hat{f}(\xi)|^2d\xi)^\frac{1}{2}.
\end{align*}

The main tool in this paper is the Littlewood-Paley decomposition of distributions into dyadic blocks of frequencies:
\begin{definition}\label{LPD}
Let $\psi(\xi)\in C^\infty_0(B_1)$ be a non-negative function so that $\psi(\xi)=1$ for $|\xi|\leq \frac{1}{2}$. Let $\varphi(\xi)$ be defined as $\varphi(\xi)=\psi(2^{-1}\xi)-\psi(\xi)$.
For given $u\in\mathcal{S}'(\mathbb{R}^3)$, the homogeneous dyadic blocks $\dot{\Delta}_k$ and the homogeneous low-frequency cut-off operator $\dot{S}_k$ are defined for all $k\in\mathbb{Z}$ by
\begin{align}\label{block}
\dot{\Delta}_ku(x)&=\mathcal{F}^{-1}(\varphi(2^{-k}\cdot)\hat{u}(\cdot))(x)=\frac{1}{(2\pi)^\frac{3}{2}}\int_{\mathbb{R}^3}e^{ix\cdot\xi}\varphi(2^{-k}\xi)\hat{u}(\xi)d\xi,\\
\dot{S}_ku(x)&=\mathcal{F}^{-1}(\psi(2^{-k}\cdot)\hat{u}(\cdot))(x)=\frac{1}{(2\pi)^\frac{3}{2}}\int_{\mathbb{R}^3}e^{ix\cdot\xi}\psi(2^{-k}\xi)\hat{u}(\xi)d\xi.
\end{align}
\end{definition}

 Throughout this paper, we will use the notation that $\tilde{\dot{\Delta}}_lu=\sum_{|l'-l|\leq2}\dot{\Delta}_{l'}u$. In the homogeneous case, the following Littlewood-Paley decomposition makes sense
\begin{align*}
u(x)=\sum_{k\in\mathbb{Z}}\dot{\Delta}_ku(x)\quad for\quad u\in\mathcal{S}_h'(\mathbb{R}^3),
\end{align*} where $\mathcal{S}_h'(\mathbb{R}^3)$ is given by
\begin{align*}
\mathcal{S}_h'(\mathbb{R}^3)=\{u\in\mathcal{S}'(\mathbb{R}^3):\lim_{k\to-\infty}\dot{S}_ku=0\}.
\end{align*}Moreover, it holds that
\begin{align*}
\dot{S}_ku(x)=\sum_{l\leq k-1}\dot{\Delta}_lu(x).
\end{align*} 

A great advantage of the localized techniques in frequency is the so-called Bernstein inequalities which will be used in the sequel.
\begin{lemma}\cite{B-C-D}
Let $\mathcal{C}$ be an annulus and $B$ a ball. A constant C exists such that
for any nonnegative integer $k$, any couple $(p,q)$ in $[1,\infty]^2$ with $1\leq p\leq q$, and any function $u\in L^p$, we have
\begin{align}
supp\hat{u}\subset\lambda B&\Rightarrow ||\nabla^k u||_{L^q}\leq C^{k+1}\lambda^{k+\frac{3}{p}-\frac{3}{q}}||u||_{L^p},\label{Ber1}\\
supp\hat{u}\subset\lambda \mathcal{C}&\Rightarrow C^{-k-1}\lambda^k||u||_{L^p}\leq||\nabla^ku||_{L^p}\leq C^{k+1}\lambda^k||u||_{L^p}\label{Ber2}
\end{align}
\end{lemma}
Based on the Littlewood-Paley decomposition and Bernstein inequality \eqref{Ber2}, the homogeneous Sobolev norm $||\cdot||_{\dot{H}^s}$ can be equivalently written as
\begin{align*}
||f||_{\dot{H}^s(\mathbb{R}^3)}=(\sum_{k\in\mathbb{Z}}2^{2ks}||\dot{\Delta}_k f||^2_{L^2})^\frac{1}{2}.
\end{align*}
When we investigate the stationary fractional Navier-Stokes equations, we need the following estimate.
\begin{lemma}\label{bds1}
Assume $0< s<1$ and $f,g\in\dot{H}^s(\mathbb{R}^3)$. Then it holds that
\begin{align}
||f g||_{\dot{H}^{2s-\frac{3}{2}}(\mathbb{R}^3)}\leq C||f||_{\dot{H}^s(\mathbb{R}^3)}||g||_{\dot{H}^s(\mathbb{R}^3)}.
\end{align}
\end{lemma}
\emph{\bf Proof.} The aim is to estimate $\{2^{k(2s-\frac{3}{2})}||\dot{\Delta}_k(fg)||_{L^2}\}_{{l}^2}$.
By Bony decomposition, it follows that
\begin{align}\label{BD}
&2^{k(2s-\frac{3}{2})}||\dot{\Delta}_k(fg)||_{L^2}\\\leq&
2^{k(2s-\frac{3}{2})}||\dot{\Delta}_k\sum_{i}\dot{\Delta}_i f\dot{S}_{i-2}g||_{L^2}+2^{k(2s-\frac{3}{2})}||\dot{\Delta}_k\sum_{i}\dot{\Delta}_i g\dot{S}_{i-2}f||_{L^2},\notag\\&+2^{k(2s-\frac{3}{2})}||\dot{\Delta}_k\sum_{i}\dot{\Delta}_i f\tilde{\dot{{\Delta}}}_{i}g||_{L^2},\notag\\
=&L_k+M_k+N_k.\notag
\end{align} Observing that
\begin{align*}
supp \mathcal{F}(\dot{\Delta}_if\dot{S}_{i-2}g)\subset\{\xi:2^{i-2}\leq|\xi|<\frac{9}{8}2^{i+1}\},
\end{align*} we obtain that
\begin{align*}
L_k\leq&2^{k(2s-\frac{3}{2})}\sum_{|i-k|\leq 2}||\dot{\Delta}_if||_{L^2}||\dot{S}_{i-2}g||_{L^\infty}\\
\leq&2^{k(2s-\frac{3}{2})}\sum_{|i-k|\leq 2}||\dot{\Delta}_if||_{L^2}\sum_{j\leq i-3}2^{\frac{3}{2}j}||\dot{\Delta}_jg||_{L^2},\notag\\
\leq&2^{(k-i)(2s-\frac{3}{2})}\sum_{|i-k|\leq 2}2^{is}||\dot{\Delta}_if||_{L^2}\sum_{j\leq i-3}2^{(\frac{3}{2}-s)(j-i)}2^{js}||\dot{\Delta}_jg||_{L^2},\notag\\
\leq& C\sum_{|i-k|\leq 2}2^{is}||\dot{\Delta}_if||_{L^2}\sum_{j\leq i-3}2^{(\frac{3}{2}-s)(j-i)}2^{js}||\dot{\Delta}_jg||_{L^2}\notag
\end{align*}here we have used Bernstein inequality \eqref{Ber1}. By using H\"older inequality and Young inequality, we get that
\begin{align}\label{Lk}
(\sum_{k}L_k^2)^\frac{1}{2}\leq& C(\sum_{i}(2^{is}||\dot{\Delta}_if||_{L^2})^4)^\frac{1}{4}\times(\sum_{j}(2^{js}||\dot{\Delta}_jg||_{L^2})^4)^\frac{1}{4},\\
\leq& C(\sum_{i}(2^{is}||\dot{\Delta}_if||_{L^2})^2)^\frac{1}{2}\times(\sum_{j}(2^{js}||\dot{\Delta}_jg||_{L^2})^2)^\frac{1}{2},\notag\\
\leq& C||f||_{\dot{H}^s}||g||_{\dot{H}^s}\notag
\end{align}here we have used the fact $l^2\hookrightarrow l^4$.
For $M_k$, similarly to the above one can obtain
\begin{align}\label{Mk}
(\sum_{k}M_k^2)^\frac{1}{2}\leq& C||f||_{\dot{H}^s}||g||_{\dot{H}^s}.
\end{align}

Next, we consider $N_k$. By using Bernstein inequality \eqref{Ber1}, we can obtain from the fact
$supp\mathcal{F}(\dot{\Delta}_if\tilde{\dot{\Delta}}_{i}g)\subset\{\xi:|\xi|\leq5\times2^{i+1}\}$ that
\begin{align*}
N_k\leq C 2^{2ks}||\dot{\Delta}_k\sum_{i\geq k-4}\dot{\Delta}_i f\tilde{\dot{{\Delta}}}_{i}g||_{L^1}
\leq C \sum_{i\geq k-4}2^{2(k-i)s}2^{is}||\dot{\Delta}_i f||_{L^2}2^{is}||\tilde{\dot{{\Delta}}}_{i}g||_{L^2}.
\end{align*}Note $s>0$, by using Young inequality, it follows that
\begin{align}\label{Nk}
(\sum_{k\in\mathbb{Z}}N_k^2)^\frac{1}{2}\leq C\sum_{i\in\mathbb{Z}}2^{is}||\dot{\Delta}_i f||_{L^2}2^{is}||\tilde{\dot{{\Delta}}}_{i}g||_{L^2}\leq C||f||_{\dot{H}^s}||g||_{\dot{H}^s}.
\end{align}Combining \eqref{BD}-\eqref{Nk} completes the proof of Lemma \ref{bds1}.
\qquad $\hfill\Box$

It is also worth pointing out that throughout this paper, we will use the following well-known facts:
\begin{align}\label{+10}
||\dot{S}_kf||_{L^p}\leq C||f||_{L^p},\quad ||f^k||_{L^p}\leq C||f||_{L^p}
\end{align}here $f^k=f-\dot{S}_kf$.

\section{Proofs of results}
\subsection{Classical stationary Navier-Stokes equations}
\emph{\bf Proof of Theorem \ref{main}.}
We define $u^k(x)$ which is the high-frequency part of $u$ as follows
\begin{align}\label{hp}
u^k(x)\doteq\sum_{l\geq k}\dot{\Delta}_lu(x)=u(x)-\dot{S}_ku(x).
\end{align}For $u^k(x)$, we have that
\begin{align}
||\nabla u^k||_{L^2}&\leq D(u)^\frac{1}{2},\label{es1}\\
||u^k||_{L^3}\leq\sum_{l\geq k}||\dot{\Delta}_lu||_{L^3}\leq\sum_{l\geq k} 2^{\frac{l}{2}}||\dot{\Delta}_lu||_{L^2}&\leq\sum_{l\geq k}2^{\frac{-l}{2}}||\nabla\dot{\Delta}_lu||_{L^2}\leq 2^{\frac{-k}{2}}D(u)^\frac{1}{2},\label{es2}
\end{align} here we have used the Bernstein inequalities \eqref{Ber1}-\eqref{Ber2} and H\"{o}lder inequality.

Next, we deduce the $L^\frac{3}{2}$ bound for $\nabla P$. Take the divergence operator on \eqref{NS} implies
\begin{align*}
-\Delta P=\nabla\cdot(u\cdot\nabla u).
\end{align*}It follows that $\nabla P=\nabla(-\Delta)^{-1}\nabla\cdot(u\cdot\nabla u)$. The boundness of Riesz transforms on $L^\frac{3}{2}$ implies
\begin{align}\label{es3}
||\nabla P||_{L^\frac{3}{2}}\leq C||u\cdot\nabla u||_{L^\frac{3}{2}}\leq C||u||_{L^6}||\nabla u||_{L^2}\leq CD(u),
\end{align}here we have used the fact $||u||_{L^6(\mathbb{R}^3)}\leq C||\nabla u||_{L^2(\mathbb{R}^3)}$.

In view of the estimates \eqref{es1}-\eqref{es3}, we can take $L^2$
inner product to \eqref{NS} with $u^k$, then obtain that
\begin{align}\label{ap}
\int_{\mathbb{R}^3}|\nabla u^k|^2dx=-\int_{\mathbb{R}^3}u\cdot\nabla u\cdot u^kdx-\int_{\mathbb{R}^3}\nabla\dot{S}_ku\cdot\nabla u^kdx.
\end{align} It is clear that
\begin{align*}
\lim_{k\to-\infty}\int_{\mathbb{R}^3}|\nabla u^k|^2dx=&\int_{\mathbb{R}^3}|\nabla u|^2dx,\\
|\lim_{k\to-\infty}\int_{\mathbb{R}^3}\nabla\dot{S}_ku\cdot\nabla u^kdx|\leq& \lim_{k\to-\infty}||\nabla\dot{S}_ku||_{L^2}||\nabla \dot{\Delta}_ku||_{L^2}=0.
\end{align*} Hence we deduce that
\begin{align}\label{ob}
\int_{\mathbb{R}^3}|\nabla u|^2dx=-\liminf_{k\to-\infty}\int_{\mathbb{R}^3}u\cdot\nabla u\cdot u^kdx.
\end{align}

We will localize the integral in \eqref{ob} into the region near the origin in frequency space. We split this integral as follows
\begin{align}\label{dec}
&\int_{\mathbb{R}^3}u\cdot\nabla u\cdot u^kdx\\
=&\int_{\mathbb{R}^3}\dot{S}_ku\cdot\nabla\dot{S}_ku\cdot u^kdx+\int_{\mathbb{R}^3}u^k\cdot\nabla\dot{S}_ku\cdot u^kdx+\int_{\mathbb{R}^3}u\cdot\nabla u^k\cdot u^kdx\notag\\
=& I_1+I_2+I_3.\notag
\end{align}Using the fact $\nabla\cdot u=0$, it is easy to see that
\begin{align}\label{van}
I_3=-\int_{\mathbb{R}^3}\nabla\cdot u\frac{|u^k|^2}{2}dx=0.
\end{align}

Next, we investigate the term $I_1$. Throughout, we will use Einstein summation convention (summing over repeated indices). By Bony decomposition, we get
\begin{align}\label{I1}
I_1=&\int_{\mathbb{R}^3}\partial_j\dot{S}_ku_i\dot{S}_ku_j u^k_idx,\\
=&\int_{\mathbb{R}^3}\partial_j\dot{S}_ku_i(\sum_{l\in\mathbb{Z}}\dot{\Delta}_l\dot{S}_ku_j\dot{S}_{l-2}u^k_i+\sum_{l\in\mathbb{Z}}\dot{\Delta}_lu^k_i\dot{S}_{l-2}\dot{S}_ku_j
+\sum_{l\in\mathbb{Z}}\dot{\Delta}_l\dot{S}_ku_j\tilde{\dot{\Delta}}_lu^k_i)dx,\notag\\
=&I_{11}+I_{12}+I_{13},\notag
\end{align} where we have used the notation $\tilde{\dot{\Delta}}_lf=\sum_{|l'-l|\leq 2}\dot{\Delta}_{l'}f$.

For $I_{11}$, we observe that $\dot{\Delta}_l\dot{S}_ku_j\dot{S}_{l-2}u^k_i\neq0$ means $l\leq k$ and $l\geq k+2$. This means that $\dot{\Delta}_l\dot{S}_ku_j\dot{S}_{l-2}u^k_i=0$ for all $l\in\mathbb{Z}$ and then
\begin{align}\label{I11}
I_{11}=\int_{\mathbb{R}^3}\partial_j\dot{S}_ku_i\sum_{l\in\mathbb{Z}}\dot{\Delta}_l\dot{S}_ku_j\dot{S}_{l-2}u^k_idx=0.
\end{align}For $I_{12}$, we first observe that $\dot{\Delta}_lu^k_i\dot{S}_{l-2}\dot{S}_ku_j\neq0$ means $l\geq k-1$ and
\begin{align}\label{es4}
supp \mathcal{F}(\dot{\Delta}_lu^k_i\dot{S}_{l-2}\dot{S}_ku_j)\subset\{\xi:2^{l-2}\leq|\xi|<\frac{9}{8}2^{l+1}\}.
\end{align} On the other hand, we also have
\begin{align}\label{es5}
supp \mathcal{F}(\partial_j\dot{S}_ku_i)\subset\{\xi:|\xi|< 2^k\}.
\end{align}Combine \eqref{es4} and \eqref{es5} implies that
\begin{align}\label{I12}
I_{12}=&\int_{\mathbb{R}^3}\partial_j\dot{S}_ku_i\sum_{l\in\mathbb{Z}}\dot{\Delta}_lu^k_i\dot{S}_{l-2}\dot{S}_ku_jdx
=\int_{\mathbb{R}^3}\partial_j\dot{S}_ku_i\sum_{l=k-1}^{k+2}\dot{\Delta}_lu^k_i\dot{S}_{l-2}\dot{S}_{k}u_jdx\\
=&\int_{\mathbb{R}^3}\sum_{l=k-1}^{k+2}\sum_{l'=l-2}^{k-1}\partial_j\dot{\Delta}_{l'}u_i\dot{\Delta}_lu^k_i\dot{S}_{l-2}\dot{S}_{k}u_jdx\notag\\=&\int_{\mathbb{R}^3}
\sum_{l=k-1}^{k+2}\sum_{l'=l-2}^{k-1}\dot{S}_{l-2}\dot{S}_{k}u\cdot\nabla\dot{\Delta}_{l'}u\cdot\dot{\Delta}_lu^kdx.\notag
\end{align} We now investigate $I_{13}$. It is not difficult to see that $\dot{\Delta}_l\dot{S}_ku_j\neq0\Rightarrow l\leq k$ and $\tilde{\dot{\Delta}}_lu^k_i\neq0\Rightarrow l\geq k-3$.
Hence it follows from the above facts that
\begin{align*}
\sum_{l\in\mathbb{Z}}\dot{\Delta}_l\dot{S}_ku_j\tilde{\dot{\Delta}}_lu^k_i=\sum_{l=k-3}^{k}\dot{\Delta}_l\dot{S}_ku_j\tilde{\dot{\Delta}}_lu^k_i.
\end{align*}From the above identity, we deduce that
\begin{align}\label{I13}
I_{13}=&\int_{\mathbb{R}^3}\partial_j\dot{S}_ku_i\sum_{l\in\mathbb{Z}}\dot{\Delta}_l\dot{S}_ku_j\tilde{\dot{\Delta}}_lu^k_i
=\int_{\mathbb{R}^3}\partial_j\dot{S}_ku_i\sum_{l=k-3}^{k}\dot{\Delta}_l\dot{S}_ku_j\tilde{\dot{\Delta}}_lu^k_idx\\
=&\int_{\mathbb{R}^3}\sum_{l=k-3}^{k}\dot{\Delta}_l\dot{S}_ku\cdot\nabla\dot{S}_ku\cdot\tilde{\dot{\Delta}}_lu^kdx.\notag
\end{align}Based on \eqref{I1}, \eqref{I11}, \eqref{I12} and \eqref{I13}, it follows that
\begin{align}\label{I1m}
I_1=\int_{\mathbb{R}^3}
\sum_{l=k-1}^{k+2}\sum_{l'=l-2}^{k-1}\dot{S}_{l-2}\dot{S}_{k}u\cdot\nabla\dot{\Delta}_{l'}u\cdot\dot{\Delta}_lu^kdx+\int_{\mathbb{R}^3}\sum_{l=k-3}^{k}\dot{\Delta}_l\dot{S}_ku\cdot\nabla\dot{S}_ku\cdot\tilde{\dot{\Delta}}_lu^kdx.
\end{align}In view of \eqref{I1m}, we see that $I_1$ is determined by the behavior of $u$ near the origin in frequency space as $k\to -\infty$.

We now restrict attention to the term $I_2$. Similarly,  by using Bony decomposition, it follows that
\begin{align}\label{I2}
I_2=&\int_{\mathbb{R}^3}\partial_j\dot{S}_ku_iu^k_ju^k_idx,\\
=&\int_{\mathbb{R}^3}\partial_j\dot{S}_ku_i(\sum_{l\in\mathbb{Z}}\dot{\Delta}_lu^k_j\dot{S}_{l-2}u^k_i+\sum_{l\in\mathbb{Z}}\dot{\Delta}_lu^k_i\dot{S}_{l-2}u^k_j
+\sum_{l\in\mathbb{Z}}\dot{\Delta}_lu^k_i\tilde{\dot{\Delta}}_{l}u^k_j)dx,\notag\\
=&I_{21}+I_{22}+I_{23}.\notag
\end{align}
For $I_{21}$, we first observe that $\dot{\Delta}_lu^k_j\dot{S}_{l-2}u^k_i\neq0\Rightarrow l\geq k+2$ and then
\begin{align}\label{es6}
supp\mathcal{F}(\dot{\Delta}_lu^k_j\dot{S}_{l-2}u^k_i)\subset\{\xi:2^{l-2}\leq|\xi|<\frac{9}{8}2^{l+1}\} \quad for\quad l\geq k+2.
\end{align}Combining \eqref{es6} and \eqref{es5} shows that
\begin{align}\label{I21}
I_{21}=&\int_{\mathbb{R}^3}\partial_j\dot{S}_ku_i\sum_{l\in\mathbb{Z}}\dot{\Delta}_lu^k_j\dot{S}_{l-2}u^k_idx\\
=&\int_{\mathbb{R}^3}\partial_j\dot{S}_ku_i\sum_{l\geq k+2}\dot{\Delta}_lu^k_j\dot{S}_{l-2}u^k_idx\notag\\
=&\int_{\mathbb{R}^3}\mathcal{F}(\partial_j\dot{S}_ku_i)\mathcal{F}(\sum_{l\geq k+2}\dot{\Delta}_lu^k_j\dot{S}_{l-2}u^k_i)d\xi\notag\\
=&0.\notag
\end{align}Applying the same arguments, we also obtain
\begin{align}\label{I22}
I_{22}=\int_{\mathbb{R}^3}\partial_j\dot{S}_ku_i\sum_{l\in\mathbb{Z}}\dot{\Delta}_lu^k_i\dot{S}_{l-2}u^k_jdx=0.
\end{align}

Finally, we consider the most complicated quantity $I_{23}$. It is not difficult to see that $\dot{\Delta}_lu^k_i\neq0\Rightarrow l\geq k-1$, we thus get that
\begin{align*}
I_{23}=\int_{\mathbb{R}^3}\partial_j\dot{S}_ku_i\sum_{l\in\mathbb{Z}}\dot{\Delta}_lu^k_i\tilde{\dot{\Delta}}_{l}u^k_jdx
=\int_{\mathbb{R}^3}\partial_j\dot{S}_ku_i\sum_{l\geq k-1}\dot{\Delta}_lu^k_i\tilde{\dot{\Delta}}_{l}u^k_jdx.
\end{align*} Since $supp\mathcal{F}(\dot{\Delta}_lu^k_i\tilde{\dot{\Delta}}_{l}u^k_j)\subset\{\xi:|\xi|<5\times2^{l+1}\}$, we can not localize $I_{23}$ into the region near origin in frequency space directly.
The key observation is that we can decompose $I_{23}$ into two parts, the first one is a low frequency part compared to $k$, the other is a high frequency part compared to $k$ and the high frequency part will vanish as $k\to-\infty$. Let $\theta\in (0,1)$ whose value will be determined later. We decompose $I_{23}$ as follows
\begin{align}\label{deI23}
I_{23}=&\int_{\mathbb{R}^3}\partial_j\dot{S}_ku_i\sum_{k-1\leq l\leq [\theta k]}\dot{\Delta}_lu^k_i\tilde{\dot{\Delta}}_{l}u^k_jdx+\int_{\mathbb{R}^3}\partial_j\dot{S}_ku_i\sum_{l\geq [\theta k]+1}\dot{\Delta}_lu^k_i\tilde{\dot{\Delta}}_{l}u^k_jdx\\
=&I_{231}+I_{232}.\notag
\end{align} The quantity $I_{232}$ can be bounded as following:
\begin{align}\label{bI232}
I_{232}&\leq||\partial_j\dot{S}_ku_i||_{L^\infty}\sum_{l\geq [\theta k]+1}||\dot{\Delta}_lu^k_i||_{L^2}||\tilde{\dot{\Delta}}_{l}u^k_j||_{L^2},\\
&\leq 2^{\frac{3}{2}k}||\dot{S}_ku_i||_{L^6}\sum_{l\geq [\theta k]+1}2^{-2l}||\nabla\dot{\Delta}_lu^k_i||_{L^2}||\nabla\tilde{\dot{\Delta}}_{l}u^k_j||_{L^2},\notag\\
&\leq C2^{(\frac{3}{2}-2\theta )k}||\nabla\dot{S}_ku_i||_{L^2}D(u)\notag
\end{align} where we have used \eqref{Ber1}-\eqref{Ber2} and the fact $\dot{H}^1(\mathbb{R}^3)\hookrightarrow L^6(\mathbb{R}^3)$. From \eqref{bI232}, We deduce
\begin{align}\label{vn}
\lim_{k\to-\infty}I_{232}=0 \quad if \quad\theta\leq\frac{3}{4}.
\end{align}Based on \eqref{dec}, \eqref{I1m}, \eqref{I21}, \eqref{I22}, \eqref{deI23} and \eqref{vn}, we obtain choosing $\theta=\frac{1}{2}$ that
\begin{align}\label{lim}
\liminf_{k\to-\infty}&\int_{\mathbb{R}^3}u\cdot\nabla u\cdot u^kdx\\=\liminf_{k\to-\infty}&\{\int_{\mathbb{R}^3}
\sum_{l=k-1}^{k+2}\sum_{l'=l-2}^{k-1}\dot{S}_{l-2}\dot{S}_{k}u\cdot\nabla\dot{\Delta}_{l'}u\cdot\dot{\Delta}_lu^kdx+\int_{\mathbb{R}^3}\sum_{l=k-3}^{k}\dot{\Delta}_l\dot{S}_ku\cdot\nabla\dot{S}_ku\cdot\tilde{\dot{\Delta}}_lu^kdx\notag\\&+\int_{\mathbb{R}^3}\sum_{l=k-1}^{[\frac{k}{2}]}\dot{\Delta}_l u^k\cdot\nabla \dot{S}_ku\cdot\tilde{\dot{\Delta}}_lu^kdx\}.\notag
\end{align}Substitute \eqref{lim} into \eqref{ob}, we then complete the proof of Theorem \ref{main}.
\qquad $\hfill\Box$

Next, we will show a Liouville type theorem based on \eqref{snc}.

\emph{\bf{Proof of Corollary \ref{coro}.}} The aim is to bound the right hand side terms in \eqref{snc} based on Bernstein inequalities \eqref{Ber1}-\eqref{Ber2} and \eqref{+10}. For the first term and the second term, we estimate them as follows.
\begin{align}\label{bsnc1}
&\int_{\mathbb{R}^3}
\sum_{l=k-1}^{k+2}\sum_{l'=l-2}^{k-1}\dot{S}_{l-2}\dot{S}_{k}u\cdot\nabla\dot{\Delta}_{l'}u\cdot\dot{\Delta}_lu^kdx+\int_{\mathbb{R}^3}\sum_{l=k-3}^{k}
\dot{\Delta}_l\dot{S}_ku\cdot\nabla\dot{S}_ku\cdot\tilde{\dot{\Delta}}_lu^kdx\\
\leq&\sum_{l=k-1}^{k+2}\sum_{l'=l-2}^{k-1}||\dot{S}_{l-2}\dot{S}_ku||_{L^\infty}||\nabla\dot{\Delta}_{l'}u||_{L^2}||\dot{\Delta}_lu^k||_{L^2}
+\sum_{l=k-3}^{k}||\nabla\dot{S}_{k}u||_{L^\infty}||\dot{\Delta}_{l}\dot{S}_ku||_{L^2}||\tilde{\dot{\Delta}}_lu^k||_{L^2}\notag\\
\leq&\sum_{l=k-1}^{k+2}\sum_{l'=l-2}^{k-1}||\dot{S}_ku||_{L^\infty}||\nabla\dot{\Delta}_{l'}u||_{L^2}||\dot{\Delta}_lu||_{L^2}
+\sum_{l=k-3}^{k}2^k||\dot{S}_{k}u||_{L^\infty}||\dot{\Delta}_{l}u||_{L^2}||\tilde{\dot{\Delta}}_lu||_{L^2}\notag\\
\leq&2^{-k}||\dot{S}_ku||_{L^\infty}(\sum_{l=k-1}^{k+2}\sum_{l'=l-2}^{k-1}2^{k-l}||\nabla\dot{\Delta}_{l'}u||_{L^2}||\nabla\dot{\Delta}_lu||_{L^2}
+\sum_{l=k-3}^{k}2^{2k-2l}||\nabla\dot{\Delta}_{l}u||_{L^2}||\nabla\tilde{\dot{\Delta}}_lu||_{L^2})\notag\\
\leq&C2^{-k}||\dot{S}_{k}u||_{L^\infty}(\sum_{l=k-1}^{k+2}\sum_{l'=l-2}^{k-1}||\nabla\dot{\Delta}_{l'}u||_{L^2}||\nabla\dot{\Delta}_lu||_{L^2}
+\sum_{l=k-3}^{k}||\nabla\dot{\Delta}_{l}u||_{L^2}||\nabla\tilde{\dot{\Delta}}_lu||_{L^2})\notag\\
\leq& C2^{-k}||\dot{S}_{k}u||_{L^\infty}||\nabla\dot{S}_{k+3}u||^2_{L^2}.\notag
\end{align}Then, we control the third term as follows.
\begin{align}\label{bsnc2}
&\int_{\mathbb{R}^3}\sum_{l=k-1}^{[\frac{k}{2}]}\dot{\Delta}_l u^k\cdot\nabla \dot{S}_ku\cdot\tilde{\dot{\Delta}}_lu^kdx,\\
\leq&||\nabla\dot{S}_ku||_{L^\infty}\sum_{l=k-1}^{[\frac{k}{2}]}||\dot{\Delta}_l u||_{L^2}||\tilde{\dot{\Delta}}_lu||_{L^2},\notag\\
\leq& C||\nabla\dot{S}_ku||_{L^\infty}\sum_{l=k-1}^{[\frac{k}{2}]}2^{-2l}||\nabla\dot{\Delta}_l u||_{L^2}||\nabla\tilde{\dot{\Delta}}_lu||_{L^2},\notag\\
\leq& C2^{-k}||\dot{S}_ku||_{L^\infty}\sum_{l=k-1}^{[\frac{k}{2}]}||\nabla\dot{\Delta}_l u||_{L^2}||\nabla\tilde{\dot{\Delta}}_lu||_{L^2},\notag\\
\leq&C2^{-k}||\dot{S}_ku||_{L^\infty}||\nabla\dot{S}_{[\frac{k}{2}]+3}u||^2_{L^2}.\notag
\end{align} Notice that $||\nabla u||^2_{L^2}=D(u)<\infty$, one can see that
\begin{align}\label{sm}
\lim_{k\to-\infty}(||\nabla\dot{S}_{k+3}u||^2_{L^2}+||\nabla\dot{S}_{[\frac{k}{2}]+1}u||^2_{L^2})=0.
\end{align}Substituting the estimates \eqref{bsnc1}-\eqref{sm} into \eqref{snc}, we deduce that if \eqref{sc1} holds then we can get that
\begin{align*}
&\int_{\mathbb{R}^3}|\nabla u|^2dx\\=&-\liminf_{k\to-\infty}\{\int_{\mathbb{R}^3}
\sum_{l=k-1}^{k+2}\sum_{l'=l-2}^{k-1}\dot{S}_{l-2}\dot{S}_{k}u\cdot\nabla\dot{\Delta}_{l'}u\cdot\dot{\Delta}_lu^kdx+\int_{\mathbb{R}^3}\sum_{l=k-3 }^{k}\dot{\Delta}_l\dot{S}_ku\cdot\nabla\dot{S}_ku\cdot\tilde{\dot{\Delta}}_lu^kdx\notag\\&+\int_{\mathbb{R}^3}\sum_{k-1\leq l\leq[\frac{k}{2}]}\dot{\Delta}_l u\cdot\nabla \dot{S}_ku\cdot\tilde{\dot{\Delta}}_ludx\}\notag\\
\leq&\liminf_{k\to-\infty} C2^{-k}||\dot{S}_ku||_{L^\infty}\times\lim_{k\to-\infty}(||\nabla\dot{S}_{k+3}u||^2_{L^2}+||\nabla\dot{S}_{[\frac{k}{2}]+1}u||^2_{L^2})\notag\\
=&0.\notag
\end{align*} It follows that $||u||_{L^6(\mathbb{R}^3)}\leq C||\nabla u||_{L^2(\mathbb{R}^3)}=0$. We thus conclude $u\equiv0$.

Furthermore, notice that
\begin{align}\label{es7}
&2^{-k}||\dot{S_k}u||_{L^\infty}\leq2^{-k}\sum_{l\leq {k-1}}||\dot{\Delta}_lu||_{L^\infty}\\
\leq& 2^{-k}\sum_{l\leq {k-1}}2^l||\dot{S}_ku||_{\dot{B}^{-1}_{\infty,\infty}}\leq C||\dot{S}_ku||_{\dot{B}^{-1}_{\infty,\infty}}\notag
\end{align}and
\begin{align}\label{es8}
&2^{-k}||\dot{S_k}u||_{L^\infty}\leq2^{-k}||\hat{u}||_{L^1(B_{2^k})}\\
\leq& 2^{k(-1+3(1-\frac{1}{r}))}||\hat{u}||_{L^r(B_{2^k})}\leq ||\hat{u}||_{L^r(B_{2^k})}\quad for \quad r\geq\frac{3}{2}.\notag
\end{align}From \eqref{es7} and \eqref{es8}, it is easy to see that \eqref{sc2} or \eqref{sc3} implies \eqref{sc1}. We thus complete the proof of Corollary \ref{coro}.
\qquad $\hfill\Box$
\subsection{Stationary fractional Navier-Stokes equations}
\emph{\bf Proof of Theorem \ref{main1}.} By Lemma \ref{bds1}, we have that
\begin{align}\label{es9}
||u_iu_j||_{\dot{H}^{2s-\frac{3}{2}}}\leq ||u||^2_{\dot{H}^s}.
\end{align} Furthermore, it also holds
\begin{align}\label{es10}
||\partial_ju^k_i||_{\dot{H}^{\frac{3}{2}-2s}}=&(\sum_{l\geq k}(2^{(\frac{3}{2}-2s)l}||\dot{\Delta}_l\partial_ju_i||_{L^2})^2)^\frac{1}{2}\\
=&(\sum_{l\geq k}(2^{(\frac{5}{2}-3s)l}2^{(s-1)l}||\dot{\Delta}_l\partial_ju_i||_{L^2})^2)^\frac{1}{2}\notag\\
\leq& C(k)(\sum_{l\geq k}(2^{(s-1)l}||\dot{\Delta}_l\partial_ju_i||_{L^2})^2)^\frac{1}{2}\notag\\
\leq& C(k)||\partial_ju_i||_{\dot{H}^{s-1}}\leq C(k)||u||_{\dot{H}^s},\notag
\end{align} here we have used the condition $s\geq\frac{5}{6}$. From \eqref{es9} and \eqref{es10}, one can see that the integral $\int_{\mathbb{R}^3}u_ju_i\partial_ju^k_idx$ is well defined. Indeed, it holds that
\begin{align}\label{es11}
\int_{\mathbb{R}^3}u_ju_i\partial_ju^k_idx=&\int_{\mathbb{R}^3}|\xi|^{2s-\frac{3}{2}}\mathcal {F}(u_ju_i)|\xi|^{\frac{3}{2}-2s}\mathcal {F}(\partial_ju^k_i)d\xi\\
\leq& C(k)||u_iu_j||_{\dot{H}^{2s-\frac{3}{2}}}||\partial_ju^k_i||_{\dot{H}^{\frac{3}{2}-2s}}\leq C(k)||u||^3_{\dot{H}^s}.\notag
\end{align} In order to take $L^2$ inner product to \eqref{FNS} with $u^k$, we need also to establish a $\dot{H}^{2s-\frac{3}{2}}$ bound for pressure $P$.
Since
\begin{align}\label{pre}
-\Delta P=\nabla\cdot\nabla\cdot(u\otimes u),
\end{align}we deduce
from the boundedness of the Riesz transform on $\dot{H}^{2s-\frac{3}{2}}$ that
\begin{align}\label{es12}
||P||_{\dot{H}^{2s-\frac{3}{2}}}\leq C||u^2||_{\dot{H}^{2s-\frac{3}{2}}}\leq C||u||^2_{\dot{H}^s}.
\end{align}
Based on \eqref{es9}, \eqref{es10} and \eqref{es12}, we now can take $L^2$ inner product to \eqref{FNS} with $u^k$ then obtain that
\begin{align}\label{ap1}
\int_{\mathbb{R}^3}|(-\Delta)^\frac{s}{2} u^k|^2dx=\int_{\mathbb{R}^3}u_j u_i\partial_j u^k_idx-\int_{\mathbb{R}^3}(-\Delta)^\frac{s}{2}\dot{S}_{k}u\cdot(-\Delta)^\frac{s}{2}u^kdx.
\end{align} Firstly,  it is not difficult to see that
\begin{align}
\lim_{k\to-\infty}\int_{\mathbb{R}^3}|(-\Delta)^\frac{s}{2} u^k|^2dx=\int_{\mathbb{R}^3}|(-\Delta)^\frac{s}{2} u|^2dx,\label{+6}\\
\lim_{k\to-\infty}\int_{\mathbb{R}^3}(-\Delta)^\frac{s}{2}\dot{S}_{k}u\cdot(-\Delta)^\frac{s}{2}u^kdx=0.\label{+7}
\end{align}Secondly,
we apply the same arguments of the proof of Theorem \ref{main} to compute the integral $\int_{\mathbb{R}^3}u_j u_i\partial_j u^k_idx$ then obtain that
\begin{align}\label{es13}
&-\int_{\mathbb{R}^3}u_j u_i\partial_j u^k_idx\\=&\int_{\mathbb{R}^3}\dot{S}_ku_j\partial_j\dot{S}_ku_i\cdot u^k_idx+\int_{\mathbb{R}^3}u^k_j\partial_j\dot{S}_ku_i\cdot u^k_idx\notag\\
=&\int_{\mathbb{R}^3}
\sum_{l=k-1}^{k+2}\sum_{l'=l-2}^{k-1}\dot{S}_{l-2}\dot{S}_{k}u_j\partial_j\dot{\Delta}_{l'}u_i\dot{\Delta}_lu^k_idx
+\int_{\mathbb{R}^3}\sum_{l=k-3}^{k}\dot{\Delta}_l\dot{S}_ku_j\partial_j\dot{S}_ku_i\tilde{\dot{\Delta}}_lu^k_idx\notag\\
&+\int_{\mathbb{R}^3}\partial_j\dot{S}_ku_i\sum_{l=k-1}^{[\frac{k}{2}]}\dot{\Delta}_lu^k_i\tilde{\dot{\Delta}}_{l}u^k_jdx+
\int_{\mathbb{R}^3}\partial_j\dot{S}_ku_i\sum_{l\geq[\frac{k}{2}]+1}\dot{\Delta}_lu^k_i\tilde{\dot{\Delta}}_{l}u^k_jdx.\notag
\end{align}
 From the Bernstein inequalities \eqref{Ber1}-\eqref{Ber2} and Young inequality, we estimate the last term in \eqref{es13} as follows
\begin{align}\label{es14}
&|\int_{\mathbb{R}^3}\partial_j\dot{S}_ku_i\sum_{l\geq [\frac{k}{2}]+1}\dot{\Delta}_lu^k_i\tilde{\dot{\Delta}}_{l}u^k_jdx|\\
\leq& ||\partial_j\dot{S}_ku_i||_{L^\infty}\sum_{l\geq [\frac{k}{2}]+1}2^{-2ls}2^{ls}||\dot{\Delta}_lu^k_i||_{L^2}2^{ls}||\tilde{\dot{\Delta}}_{l}u^k_j||_{L^2}\notag\\
\leq& C2^{k+\frac{3-2s}{2}k- sk}||\dot{S}_ku_i||_{L^\frac{6}{3-2s}}\sum_{l\geq [\frac{k}{2}]+1}2^{-2(l-\frac{k}{2})s}2^{ls}||\dot{\Delta}_lu^k_i||_{L^2}2^{ls}||\tilde{\dot{\Delta}}_{l}u^k_j||_{L^2}\notag\\
\leq& C2^{k(\frac{5}{2}-2s)}||\dot{S}_ku_i||_{\dot{H}^s}||u_i||_{\dot{H}^s}||u_j||_{\dot{H}^s}.\notag
\end{align}Since $s<1$, we obtain from \eqref{es14}
that
\begin{align}\label{es15}
\liminf_{k\to-\infty}\int_{\mathbb{R}^3}\partial_j\dot{S}_ku_i\sum_{l\geq [\frac{k}{2}]+1}\dot{\Delta}_lu_i\tilde{\dot{\Delta}}_{l}u_jdx=0.
\end{align}Finally, collecting \eqref{ap1}-\eqref{es15}, we conclude that
\begin{align}\label{+8}
&\int_{\mathbb{R}^3}|(-\Delta)^\frac{s}{2} u|^2dx\\
=&\liminf_{k\to-\infty}\{\int_{\mathbb{R}^3}
\sum_{l=k-1}^{k+2}\sum_{l'=l-2}^{k-1}\dot{S}_{l-2}\dot{S}_{k}u_j\partial_j\dot{\Delta}_{l'}u_i\dot{\Delta}_lu^k_idx+\int_{\mathbb{R}^3}\sum_{l=k-3 }^{k}\dot{\Delta}_l\dot{S}_ku_j\partial_j\dot{S}_ku_i\tilde{\dot{\Delta}}_lu^k_idx\notag\\
&+\int_{\mathbb{R}^3}\partial_j\dot{S}_ku_i\sum_{l=k-1}^{[\frac{k}{2}]}\dot{\Delta}_lu^k_i\tilde{\dot{\Delta}}_{l}u^k_jdx\}.\notag
\end{align} This completes
the proof of Theorem \ref{main1}.
\qquad $\hfill\Box$

\emph{\bf{Proof of Corollary \ref{coro1}.}} The aim is to bound the right hand side terms in \eqref{snc1} based on Bernstein inequalities \eqref{Ber1}-\eqref{Ber2} and \eqref{+10}. For the first term and the second term, we estimate them as follows.
\begin{align}\label{es16}
&\int_{\mathbb{R}^3}\sum_{l=k-1}^{k+2}\sum_{l'=l-2}^{k-1}\dot{S}_{l-2}\dot{S}_{k}u\cdot\nabla\dot{\Delta}_{l'}u\cdot\dot{\Delta}_lu^kdx+\int_{\mathbb{R}^3}\sum_{l=k-3}^{k}
\dot{\Delta}_l\dot{S}_ku\cdot\nabla\dot{S}_ku\cdot\tilde{\dot{\Delta}}_lu^kdx\\
\leq&\sum_{l=k-1}^{k+2}\sum_{l'=l-2}^{k-1}||\dot{S}_{l-2}\dot{S}_ku||_{L^\infty}||\nabla\dot{\Delta}_{l'}u||_{L^2}||\dot{\Delta}_lu^k||_{L^2}
+\sum_{l=k-3}^{k}||\nabla\dot{S}_{k}u||_{L^\infty}||\dot{\Delta}_{l}\dot{S}_ku||_{L^2}||\tilde{\dot{\Delta}}_lu^k||_{L^2}\notag\\
\leq&\sum_{l=k-1}^{k+2}\sum_{l'=l-2}^{k-1}2^{l'}||\dot{S}_ku||_{L^\infty}||\dot{\Delta}_{l'}u||_{L^2}||\dot{\Delta}_lu||_{L^2}
+\sum_{l=k-3}^{k}2^k||\dot{S}_{k}u||_{L^\infty}||\dot{\Delta}_{l}u||_{L^2}||\tilde{\dot{\Delta}}_lu||_{L^2}\notag\\
\leq&C2^{k(1-2s)}||\dot{S}_ku||_{L^\infty}(\sum_{l=k-1}^{k+2}\sum_{l'=l-2}^{k-1}2^{l's}||\dot{\Delta}_{l'}u||_{L^2}2^{ls}||\dot{\Delta}_lu||_{L^2}
+\sum_{l=k-3}^{k}2^{ls}||\dot{\Delta}_{l}u||_{L^2}2^{ls}||\tilde{\dot{\Delta}}_lu||_{L^2})\notag\\
\leq& C2^{k(1-2s)}||\dot{S}_{k}u||_{L^\infty}||\dot{S}_{k+3}u||^2_{\dot{H}^s}.\notag
\end{align}Then, from Bernstein inequalities \eqref{Ber1} and Young inequality, we control the third term as follows.
\begin{align}\label{es17}
&\int_{\mathbb{R}^3}\sum_{k-1\leq l\leq[\frac{k}{2}]}\dot{\Delta}_l u^k\cdot\nabla \dot{S}_ku\cdot\tilde{\dot{\Delta}}_lu^kdx,\\
\leq&||\nabla\dot{S}_ku||_{L^\infty}\sum_{k-1\leq l\leq[\frac{k}{2}]}||\dot{\Delta}_l u||_{L^2}||\tilde{\dot{\Delta}}_lu||_{L^2},\notag\\
\leq& ||\nabla\dot{S}_ku||_{L^\infty}\sum_{k-1\leq l\leq[\frac{k}{2}]}2^{-2ls}2^{ls}||\nabla\dot{\Delta}_l u||_{L^2}2^{ls}||\nabla\tilde{\dot{\Delta}}_lu||_{L^2},\notag\\
\leq& C2^{k(1-2s)}||\dot{S}_ku||_{L^\infty}\sum_{k-1\leq l\leq[\frac{k}{2}]}2^{2(k-l)s}2^{ls}||\dot{\Delta}_l u||_{L^2}2^{ls}||\tilde{\dot{\Delta}}_lu||_{L^2},\notag\\
\leq&C2^{k(1-2s)}||\dot{S}_ku||_{L^\infty}||\dot{S}_{[\frac{k}{2}]+3}u||^2_{\dot{H}^s}.\notag
\end{align} In view of \eqref{FD}, one can deduce that
\begin{align}\label{es18}
\lim_{k\to-\infty}(||\dot{S}_{k+3}u||^2_{\dot{H}^s}+||\dot{S}_{[\frac{k}{2}]+3}u||^2_{\dot{H}^s})=0.
\end{align}Substituting \eqref{es16}-\eqref{es17} into the identity \eqref{snc1}, we get
\begin{align*}
&\int_{\mathbb{R}^3}|(-\Delta)^\frac{s}{2} u|^2dx\\
=&\liminf_{k\to-\infty}\{\int_{\mathbb{R}^3}
\sum_{l=k-1}^{k+2}\sum_{l'=l-2}^{k-1}\dot{S}_{l-2}\dot{S}_{k}u\cdot\nabla\dot{\Delta}_{l'}u\cdot\dot{\Delta}_lu^kdx+\int_{\mathbb{R}^3}\sum_{l=k-3 }^{k}\dot{\Delta}_l\dot{S}_ku\cdot\nabla\dot{S}_ku\cdot\tilde{\dot{\Delta}}_lu^kdx\notag\\
&+\int_{\mathbb{R}^3}\partial_j\dot{S}_ku_i\sum_{l=k-1}^{[\frac{k}{2}]}\dot{\Delta}_lu^k_i\tilde{\dot{\Delta}}_{l}u^k_jdx\}\notag\\
\leq& C\liminf_{k\to-\infty}2^{k(1-2s)}||\dot{S}_ku||_{L^\infty}\lim_{k\to-\infty}(||\dot{S}_{k+3}u||^2_{\dot{H}^s}+||\dot{S}_{[\frac{k}{2}]+3}u||^2_{\dot{H}^s}).\notag
\end{align*} From \eqref{fsc1} and \eqref{es18}, we obtain that
\begin{align*}
&\int_{\mathbb{R}^3}|(-\Delta)^\frac{s}{2} u|^2dx\\\leq& C\liminf_{k\to-\infty}2^{k(1-2s)}||\dot{S}_ku||_{L^\infty}\lim_{k\to-\infty}(||\dot{S}_{k+3}u||^2_{\dot{H}^s}+||\dot{S}_{[\frac{k}{2}]+3}u||^2_{\dot{H}^s})\\=&0.
\end{align*}
Notice that $||u||_{L^{\frac{6}{3-2s}}(\mathbb{R}^3)}\leq C ||u||_{\dot{H}^s(\mathbb{R}^3)}=0$, we thus get $u\equiv0$.

Furthermore, notice $s\geq\frac{5}{6}>\frac{1}{2}$, it follows from Young inequality and H\"{o}lder inequality that
\begin{align}\label{es19}
&2^{k(1-2s)}||\dot{S_k}u||_{L^\infty}\leq2^{k(1-2s)}\sum_{l\leq {k-1}}||\dot{\Delta}_lu||_{L^\infty}\\
\leq& 2^{k(1-2s)}\sum_{l\leq {k-1}}2^{l(2s-1)}2^{l(1-2s)}||\dot{\Delta}_lu||_{L^\infty}\leq C||\dot{S}_ku||_{\dot{B}^{1-2s}_{\infty,\infty}}\notag
\end{align}and
\begin{align}\label{es20}
&2^{k(1-2s)}||\dot{S_k}u||_{L^\infty}\leq2^{k(1-2s)}||\hat{u}||_{L^1(B_{2^k})}\\
\leq& 2^{k(1-2s+3(1-\frac{1}{r}))}||\hat{u}||_{L^r(B_{2^k})}\leq ||\hat{u}||_{L^r(B_{2^k})}\quad for \quad r\geq\frac{3}{4-2s}.\notag
\end{align}From \eqref{es19} and \eqref{es20}, it is easy to see that \eqref{fsc2} or \eqref{fsc3} implies \eqref{fsc1}. We thus complete the proof of Corollary \ref{coro1}.
\qquad $\hfill\Box$

\emph{\bf{Proof of Corollary \ref{coro2}.}}
Firstly, we have by using \eqref{Ber1} and the embedding theorem $||u||_{L^\frac{6}{3-2s}(\mathbb{R}^3)}\leq C||u||_{\dot{H}^s(\mathbb{R}^3)}$ that
\begin{align}\label{es21}
2^{k(1-2s)}||\dot{S}_ku||_{L^\infty}\leq C 2^{k(1-2s)}2^{k\frac{3-2s}{2}}||\dot{S}_ku||_{L^\frac{6}{3-2s}}\leq C2^{k(\frac{5}{2}-3s)}||u||_{\dot{H}^s}.
\end{align}When $s=\frac{5}{6}$, it is clear that
\begin{align}\label{es22}
2^{k(1-2\frac{5}{6})}||\dot{S}_ku||_{L^\infty}\leq C||u||_{\dot{H}^\frac{5}{6}}<\infty.
\end{align} We thus deduce using Corollary \ref{coro1} that $u\equiv0$.
\qquad $\hfill\Box$

\emph{\bf{Proof of Theorem \ref{main2}.}} When $\frac{1}{2}< s<\frac{5}{6}$, we know from Lemma \ref{bds1} that $u_i u_j\in \dot{H}^{2s-\frac{3}{2}}$
and $\partial_j u_i\in \dot{H}^{s-1}$, but we can not deduce that $\partial _j u^k_i\in \dot{H}^{\frac{3}{2}-2s}$ due to $s-1<\frac{3}{2}-2s$. This means that the integral $\int_{\mathbb{R}^3}u_j u_i\partial_j u^k_idx$
may be not well-defined. However, it is clear that $\partial_j\dot{S}_ku_i\in \dot{H}^{\frac{3}{2}-2s}$, we thus apply the same procedure of the
proof of Theorem \ref{main1} to obtain by taking $L^2$ inner product to \eqref{FNS} with $\dot{S}_ku$ that
\begin{align}\label{es23}
&\int_{\mathbb{R}^3}|(-\Delta)^\frac{s}{2} \dot{S}_ku|^2dx\\
=&\int_{\mathbb{R}^3}u_j u_i\partial_j \dot{S}_ku_idx
-\int_{\mathbb{R}^3}(-\Delta)^\frac{s}{2}\dot{S}_ku\cdot(-\Delta)^\frac{s}{2}u^kdx\notag\\
=&\int_{\mathbb{R}^3}\dot{S}_ku_j\partial_j{S}_ku_i\cdot u^k_idx+\int_{\mathbb{R}^3}u^k_j\partial_j{S}_ku_i\cdot u^k_idx-\int_{\mathbb{R}^3}(-\Delta)^\frac{s}{2}\dot{S}_ku\cdot(-\Delta)^\frac{s}{2}u^kdx\notag
\end{align}where we have used the fact $\int_{\mathbb{R}^3}u_j \dot{S}_ku_i\partial_j \dot{S}_ku_idx=0$ by integrating by parts and $\nabla\cdot u=0$. Taking $k\to+\infty$, we first see that
\begin{align}\label{+9}
\lim_{k\to\infty}|\int_{\mathbb{R}^3}(-\Delta)^\frac{s}{2}\dot{S}_ku\cdot(-\Delta)^\frac{s}{2}u^kdx|\leq ||u||_{\dot{H}^s}\lim_{k\to\infty}||u^k||_{\dot{H}^s}=0.
\end{align}Secondly, observing that the first term and the second term in the second line of \eqref{es23} are the same as the terms in the second line of \eqref{es13}, we
apply the same computations of \eqref{es13} to obtain that
\begin{align}\label{es24}
&\int_{\mathbb{R}^3}|(-\Delta)^\frac{s}{2}u|^2dx=\liminf_{k\to\infty}\int_{\mathbb{R}^3}|(-\Delta)^\frac{s}{2} \dot{S}_ku|^2dx\\
=&\liminf_{k\to\infty}\{\int_{\mathbb{R}^3}
\sum_{l=k-1}^{k+2}\sum_{l'=l-2}^{k-1}\dot{S}_{l-2}\dot{S}_{k}u_j\partial_j\dot{\Delta}_{l'}u_i\dot{\Delta}_lu^k_idx
+\int_{\mathbb{R}^3}\sum_{l=k-3}^{k}\dot{\Delta}_l\dot{S}_ku_j\partial_j\dot{S}_ku_i\tilde{\dot{\Delta}}_lu^k_idx\notag\\
&+\int_{\mathbb{R}^3}\partial_j\dot{S}_ku_i\sum_{l\geq k-1}\dot{\Delta}_lu^k_i\tilde{\dot{\Delta}}_{l}u^k_jdx\}\notag\\
=& \liminf_{k\to\infty}(J_1+J_2+J_3).\notag
\end{align}Our aim is to localize the integral in \eqref{es24} into the region with high frequency. For this purpose, we need to eliminate the low-frequency parts coming from
$\dot{S}_{l-2}\dot{S}_ku_j$ and $\dot{S}_ku_i$.
 Taking $0<\theta=\frac{4s-2}{3-2s}<1$, we decompose $J_1$ as follows
\begin{align}\label{es25}
J_1=&\int_{\mathbb{R}^3}
\sum_{l=k-1}^{k+2}\sum_{l'=l-2}^{k-1}\dot{S}_{l-2}\dot{S}_{k}u_j\partial_j\dot{\Delta}_{l'}u_i\dot{\Delta}_lu^k_idx\\
=&\int_{\mathbb{R}^3}
\sum_{l=k-1}^{k+2}\sum_{l'=l-2}^{k-1}\dot{S}_{l-2}\dot{S}_{[\theta k]}u_j\partial_j\dot{\Delta}_{l'}u_i\dot{\Delta}_lu^k_idx\notag\\
&+\int_{\mathbb{R}^3}
\sum_{l=k-1}^{k+2}\sum_{l'=l-2}^{k-1}\sum_{l''=[\theta k]}^{k-1}\dot{S}_{l-2}\dot{\Delta}_{l''}u_j\partial_j\dot{\Delta}_{l'}u_i\dot{\Delta}_lu^k_idx.\notag
\end{align}From Bernstein inequality \eqref{Ber1}, \eqref{+10} and the fact $\dot{H}^s(\mathbb{R}^3)\hookrightarrow L^\frac{6}{3-2s}(\mathbb{R}^3)$, we deduce that
\begin{align*}
&\lim_{k\to\infty}|\int_{\mathbb{R}^3}
\sum_{l=k-1}^{k+2}\sum_{l'=l-2}^{k-1}\dot{S}_{l-2}\dot{S}_{[\theta k]}u_j\partial_j\dot{\Delta}_{l'}u_i\dot{\Delta}_lu^k_idx|\\
\leq & \lim_{k\to\infty}C\sum_{l=k-1}^{k+2}\sum_{l'=l-2}^{k-1}2^{l'}||\dot{S}_{[\theta k]}u_j||_{L^\infty}||\dot{\Delta}_{l'}u_i||_{L^2}||\dot{\Delta}_lu_i||_{L^2}\notag\\
\leq & \lim_{k\to\infty}C2^{k((\frac{3}{2}-s)\theta+1)}||\dot{S}_{[\theta k]}u_j||_{L^\frac{6}{3-2s}}\sum_{l=k-1}^{k+2}\sum_{l'=l-2}^{k-1}||\dot{\Delta}_{l'}u_i||_{L^2}||\dot{\Delta}_lu_i||_{L^2}\notag\\
\leq& \lim_{k\to\infty}C2^{k((\frac{3}{2}-s)\theta+1-2s)}||u_j||_{\dot{H}^s}\sum_{l=k-1}^{k+2}\sum_{l'=l-2}^{k-1}2^{l's}||\dot{\Delta}_{l'}u_i||_{L^2}2^{ls}||\dot{\Delta}_lu_i||_{L^2},\notag
\end{align*}where we have used the facts $|l'-k|\leq3$ and $|l-k|\leq2$. Then notice that $(\frac{3}{2}-s)\theta+1-2s=0$ and $\sum_{l=k-1}^{k+2}\sum_{l'=l-2}^{k-1}2^{l's}||\dot{\Delta}_{l'}u_i||_{L^2}2^{ls}||\dot{\Delta}_lu_i||_{L^2}\to 0$ as $k\to\infty$, we have
\begin{align}\label{es26}
\lim_{k\to\infty}|\int_{\mathbb{R}^3}
\sum_{l=k-1}^{k+2}\sum_{l'=l-2}^{k-1}\dot{S}_{l-2}\dot{S}_{[\theta k]}u_j\partial_j\dot{\Delta}_{l'}u_i\dot{\Delta}_lu^k_idx|=0.
\end{align}
Next, we decompose $J_2$ as follows
\begin{align}\label{es27}
J_2=&\int_{\mathbb{R}^3}\sum_{l=k-3}^{k}\dot{\Delta}_l\dot{S}_ku_j\partial_j\dot{S}_ku_i\tilde{\dot{\Delta}}_lu^k_idx\\
=&\int_{\mathbb{R}^3}\sum_{l=k-3}^{k}\dot{\Delta}_l\dot{S}_ku_j\partial_j\dot{S}_{[\frac{k}{2}]}u_i\tilde{\dot{\Delta}}_lu^k_idx+
\int_{\mathbb{R}^3}\sum_{l=k-3}^{k}\sum_{l''=[\frac{k}{2}]}^{k-1}\dot{\Delta}_l\dot{S}_ku_j\partial_j\dot{\Delta}_{l''}u_i\tilde{\dot{\Delta}}_lu^k_idx.\notag
\end{align}
From Bernstein inequality \eqref{Ber1}, \eqref{+10} and the fact $\dot{H}^s(\mathbb{R}^3)\hookrightarrow L^\frac{6}{3-2s}(\mathbb{R}^3)$, we also obtain
\begin{align*}
&\lim_{k\to\infty}|\int_{\mathbb{R}^3}\sum_{l=k-3}^{k}\dot{\Delta}_l\dot{S}_ku_j\partial_j\dot{S}_{[\frac{k}{2}]}u_i\tilde{\dot{\Delta}}_lu^k_idx|\\
\leq &\lim_{k\to\infty}C2^\frac{k}{2}||\dot{S}_{[\frac{k}{2}]}u_i||_{L^\infty}\sum_{l=k-3}^{k}||\dot{\Delta}_{l}u_j||_{L^2}||\tilde{\dot{\Delta}}_lu_i||_{L^2}\notag\\
\leq&\lim_{k\to\infty}C2^{k(\frac{5}{4}-\frac{s}{2})}||\dot{S}_{[\frac{k}{2}]}u_i||_{L^\frac{6}{3-2s}}\sum_{l=k-3}^{k}||\dot{\Delta}_{l}u_j||_{L^2}||\tilde{\dot{\Delta}}_lu_i||_{L^2}\notag\\
\leq&\lim_{k\to\infty}C2^{k(\frac{5}{4}-\frac{5}{2}s)}||u_i||_{\dot{H}^s}\sum_{l=k-3}^{k}2^{ls}||\dot{\Delta}_{l}u_j||_{L^2}2^{ls}||\tilde{\dot{\Delta}}_lu_i||_{L^2},\notag
\end{align*}where we have used the fact $|l-k|\leq3$. Notice that $\sum_{l=k-3}^{k}2^{ls}||\dot{\Delta}_{l}u_j||_{L^2}2^{ls}||\tilde{\dot{\Delta}}_lu_i||_{L^2}\to0$ as $k\to\infty$ and $s\geq\frac{1}{2}$,
we thus get
\begin{align}\label{es28}
\lim_{k\to\infty}|\int_{\mathbb{R}^3}\sum_{l=k-3}^{k}\dot{\Delta}_l\dot{S}_ku_j\partial_j\dot{S}_{[\frac{k}{2}]}u_i\tilde{\dot{\Delta}}_lu^k_idx|=0.
\end{align}
For $J_3$, applying the similar computations of $J_2$ to $J_3$, we can get that
\begin{align}\label{es29}
J_3=&\int_{\mathbb{R}^3}\partial_j\dot{S}_ku_i\sum_{l\geq k-1}\dot{\Delta}_lu^k_i\tilde{\dot{\Delta}}_{l}u^k_jdx\\
=&\int_{\mathbb{R}^3}\partial_j\dot{S}_{[\frac{k}{2}]}u_i\sum_{l\geq k-1}\dot{\Delta}_lu^k_i\tilde{\dot{\Delta}}_{l}u^k_jdx
+\sum_{m=[\frac{k}{2}]}^{{k-1}}\partial_j\dot{\Delta}_mu_i\sum_{l\geq k-1}\dot{\Delta}_lu^k_i\tilde{\dot{\Delta}}_{l}u^k_jdx.\notag
\end{align}
Applying Bernstein inequality \eqref{Ber1}, \eqref{+10} and the fact $\dot{H}^s(\mathbb{R}^3)\hookrightarrow L^\frac{6}{3-2s}(\mathbb{R}^3)$, we have that
\begin{align*}
&\lim_{k\to\infty}|\int_{\mathbb{R}^3}\partial_j\dot{S}_{[\frac{k}{2}]}u_i\sum_{l\geq k-1}\dot{\Delta}_lu^k_i\tilde{\dot{\Delta}}_{l}u^k_jdx|\\
\leq&\lim_{k\to\infty}2^\frac{k}{2}||\dot{S}_{[\frac{k}{2}]}u_i||_{L^\infty}\sum_{l\geq k-1}||\dot{\Delta}_lu_i||_{L^2}||\tilde{\dot{\Delta}}_{l}u_j||_{L^2}\notag\\
\leq&\lim_{k\to\infty}2^{k(\frac{5}{4}-\frac{s}{2})}||\dot{S}_{[\frac{k}{2}]}u_i||_{L^\frac{6}{3-2s}}\sum_{l\geq k-1}2^{-2ls}2^{ls}||\dot{\Delta}_lu_i||_{L^2}2^{ls}||\tilde{\dot{\Delta}}_{l}u_j||_{L^2}\notag\\
\leq &\lim_{k\to\infty}2^{k(\frac{5}{4}-\frac{5s}{2})}||u_i||_{\dot{H}^s}\sum_{l\geq k-1}2^{ls}||\dot{\Delta}_lu_i||_{L^2}2^{ls}||\tilde{\dot{\Delta}}_{l}u_j||_{L^2}.\notag
\end{align*}From the fact $s\geq\frac{1}{2}$ and $\lim_{k\to\infty}\sum_{l\geq k-1}2^{ls}||\dot{\Delta}_{l}u_i||_{L^2}2^{ls}||\tilde{\dot{\Delta}}_lu_j||_{L^2}=0$, it follows that
\begin{align}\label{es30}
\lim_{k\to\infty}|\int_{\mathbb{R}^3}\partial_j\dot{S}_{[\frac{k}{2}]}u_i\sum_{l\geq k-1}\dot{\Delta}_lu^k_i\tilde{\dot{\Delta}}_{l}u^k_jdx|=0.
\end{align}
Substituting \eqref{es25}-\eqref{es30} into \eqref{es24} implies that
\begin{align}\label{es31}
&\int_{\mathbb{R}^3}|(-\Delta)^\frac{s}{2}u|^2dx\\=&\liminf_{k\to\infty}\{\int_{\mathbb{R}^3}
\sum_{l=k-1}^{k+2}\sum_{l'=l-2}^{k-1}\sum_{l''=[\theta k]}^{k-1}\dot{S}_{l-2}\dot{\Delta}_{l''}u_j\partial_j\dot{\Delta}_{l'}u_i\dot{\Delta}_lu^k_idx\notag\\
&+\int_{\mathbb{R}^3}\sum_{l=k-3}^{k}\sum_{l''=[\frac{k}{2}]}^{k-1}\dot{\Delta}_l\dot{S}_ku_j\partial_j\dot{\Delta}_{l''}u_i\tilde{\dot{\Delta}}_lu^k_idx
+\int_{\mathbb{R}^3}\sum_{l''=[\frac{k}{2}]}^{k-1}\sum_{k-1\leq l}\dot{\Delta}_lu^k_i\tilde{\dot{\Delta}}_{l}u^k_j\partial_j\dot{\Delta}_{l''}u_idx\}.\notag
\end{align} This completes the proof of Theorem \ref{main2}.
\qquad $\hfill\Box$

\emph{\bf{Proof of Corollary \ref{coro3}.}}
Next, we will establish the Liouville type theorem based on the identity \eqref{snc2}. From Bernstein inequalities \eqref{Ber1}-\eqref{Ber2} and \eqref{+10}, it follows that
\begin{align}
&|\int_{\mathbb{R}^3}
\sum_{l=k-1}^{k+2}\sum_{l'=l-2}^{k-1}\sum_{l''=[\theta k]}^{k-1}\dot{S}_{l-2}\dot{\Delta}_{l''}u_j\partial_j\dot{\Delta}_{l'}u_i\dot{\Delta}_lu^k_idx|\label{es32}\\
\leq&C||\sum_{l''=[\theta k]}^{k-1}\dot{S}_{l-2}\dot{\Delta}_{l''}u_j||_{L^\infty}\sum_{l=k-1}^{k+2}\sum_{l'=l-2}^{k-1}2^{l'}||\dot{\Delta}_{l'}u_i||_{L^2}||\dot{\Delta}_lu_i||_{L^2}\notag\\
\leq& C||\sum_{l''=[\theta k]}^{k-1}\dot{S}_{l-2}\dot{\Delta}_{l''}u_j||_{L^\infty}\sum_{l=k-1}^{k+2}\sum_{l'=l-2}^{k-1}2^{l'-l's-ls}2^{l's}||\dot{\Delta}_{l'}u_i||_{L^2}2^{ls}||\dot{\Delta}_lu_i||_{L^2}\notag\\
\leq& C2^{k(1-2s)} ||\sum_{l''=[\theta k]}^{k-1}\dot{\Delta}_{l''}u_j||_{L^\infty}\sum_{l=k-1}^{k+2}\sum_{l'=l-2}^{k-1}2^{l's}||\dot{\Delta}_{l'}u_i||_{L^2}2^{ls}||\dot{\Delta}_lu_i||_{L^2},\notag
\end{align}
\begin{align}
&|\int_{\mathbb{R}^3}\sum_{l=k-3}^{k}\sum_{l''=[\frac{k}{2}]}^{k-1}\dot{\Delta}_l\dot{S}_ku_j\partial_j\dot{\Delta}_{l''}u_i\tilde{\dot{\Delta}}_lu^k_idx|\label{es33}\\
\leq&||\partial_j\sum_{l''=[\frac{k}{2}]}^{k-1}\dot{\Delta}_{l''}u_i||_{L^\infty}\sum_{l=k-3}^{k}||\dot{\Delta}_{l}u_j||_{L^2}||\tilde{\dot{\Delta}}_lu_i||_{L^2}\notag\\
\leq&C2^k||\sum_{l''=[\frac{k}{2}]}^{k-1}\dot{\Delta}_{l''}u_i||_{L^\infty}\sum_{l=k-3}^{k}||\dot{\Delta}_{l}u_j||_{L^2}||\tilde{\dot{\Delta}}_lu_i||_{L^2}\notag\\
\leq&C2^{k(1-2s)}||\sum_{l''=[\frac{k}{2}]}^{k-1}\dot{\Delta}_{l''}u_i||_{L^\infty}\sum_{l=k-3}^{k}2^{ls}||\dot{\Delta}_{l}u_j||_{L^2}2^{ls}||\tilde{\dot{\Delta}}_lu_i||_{L^2},\notag
\end{align}
\begin{align}
&|\int_{\mathbb{R}^3}\sum_{l''=[\frac{k}{2}]}^{k-1}\sum_{k-1\leq l}\dot{\Delta}_lu^k_i\tilde{\dot{\Delta}}_{l}u^k_j\partial_j\dot{\Delta}_{l''}u_idx|\label{es34}\\
\leq&||\partial_j\sum_{l''=[\frac{k}{2}]}^{k-1}\dot{\Delta}_{l''}u_i||_{L^\infty}\sum_{k-1\leq l}||\dot{\Delta}_lu_i||_{L^2}||\tilde{\dot{\Delta}}_{l}u_j||_{L^2}\notag\\
\leq& C2^{k(1-2s)}||\sum_{l''=[\frac{k}{2}]}^{k-1}\dot{\Delta}_{l''}u_i||_{L^\infty}\sum_{k-1\leq l}2^{ls}||\dot{\Delta}_lu_i||_{L^2}2^{ls}||\tilde{\dot{\Delta}}_{l}u_j||_{L^2}.\notag
\end{align} Substituting the estimates \eqref{es32}-\eqref{es34} into \eqref{es31}, we conclude that
\begin{align}\label{es36}
&\int_{\mathbb{R}^3}|(-\Delta)^\frac{s}{2}u|^2dx\\\leq&\liminf_{k\to\infty}C 2^{k(1-2s)}(||\sum_{l''=[\theta k]}^{k-1}\dot{\Delta}_{l''}u_j||_{L^\infty}+||\sum_{l''=[\frac{k}{2}]}^{k-1}\dot{\Delta}_{l''}u_i||_{L^\infty})\notag\\
&\times\lim_{k\to\infty}(\sum_{l=k-1}^{k+2}\sum_{l'=l-2}^{k-1}2^{l's}||\dot{\Delta}_{l'}u_i||_{L^2}2^{ls}||\dot{\Delta}_lu_i||_{L^2}+
\sum_{l=k-3}^{k}2^{ls}||\dot{\Delta}_{l}u_j||_{L^2}2^{ls}||\tilde{\dot{\Delta}}_lu_i||_{L^2}\notag\\
&+\sum_{k-1\leq l}2^{ls}||\dot{\Delta}_lu_i||_{L^2}2^{ls}||\tilde{\dot{\Delta}}_{l}u_j||_{L^2}).\notag
\end{align} Since $||u||^2_{\dot{H}^s}=\int_{\mathbb{R}^3}|(-\Delta)^{\frac{s}{2}}u|^2dx<\infty$, it is clear that
\begin{align}\label{es35}
&\lim_{k\to\infty}\sum_{l=k-1}^{k+2}\sum_{l'=l-2}^{k-1}2^{l's}||\dot{\Delta}_{l'}u_i||_{L^2}2^{ls}||\dot{\Delta}_lu_i||_{L^2}=0,\notag\\
&\lim_{k\to\infty}\sum_{l=k-3}^{k}2^{ls}||\dot{\Delta}_{l}u_j||_{L^2}2^{ls}||\tilde{\dot{\Delta}}_lu_i||_{L^2}=0,\\
&\lim_{k\to\infty}\sum_{k-1\leq l}2^{ls}||\dot{\Delta}_lu_i||_{L^2}2^{ls}||\tilde{\dot{\Delta}}_{l}u_j||_{L^2}=0.\notag
\end{align} Substituting \eqref{es35} into \eqref{es36}, we deduce that if \eqref{+1} is valid then it follows
\begin{align*}
\int_{\mathbb{R}^3}|(-\Delta)^{\frac{s}{2}}u|^2dx=0.
\end{align*}Consequently, the embedding theorem $||u||_{L^{\frac{6}{3-2s}}(\mathbb{R}^3)}\leq C||u||_{\dot{H}^s(\mathbb{R}^3)}$ implies $u\equiv0$.

Furthermore, from the fact $s>\frac{1}{2}$ and Young inequality, it follows that
\begin{align}\label{+4}
&2^{k(1-2s)}||\sum_{l=[\theta k]}^{k-1}\dot{\Delta}_lu_j||_{L^\infty}\\
\leq& 2^{k(1-2s)}\sum_{l=[\theta k]}^{k-1}2^{l(2s-1)}2^{l(1-2s)}||\dot{\Delta}_lu_j||_{L^\infty}\notag\\
\leq& C2^{k(1-2s)}\sum_{l=[\theta k]}^{k-1}2^{l(2s-1)}||u^{[\theta k]}||_{\dot{B}^{1-2s}_{\infty,\infty}}\notag\\
\leq & C||u^{[\theta k]}||_{\dot{B}^{1-2s}_{\infty,\infty}}\notag
\end{align}and
\begin{align}\label{+5}
&2^{k(1-2s)}||\sum_{l=[\frac{k}{2}]}^{k-1}\dot{\Delta}_lu_i||_{L^\infty}\\
\leq& 2^{k(1-2s)}\sum_{l=[\frac{k}{2}]}^{k-1}2^{l(2s-1)}2^{l(1-2s)}||\dot{\Delta}_lu_j||_{L^\infty}\notag\\
\leq& 2^{k(1-2s)}\sum_{l=[\frac{k}{2}]}^{k-1}2^{l(2s-1)}||u^{[\frac{k}{2}]}||_{\dot{B}^{1-2s}_{\infty,\infty}}\notag\\
\leq & C||u^{[\frac{k}{2}]}||_{\dot{B}^{1-2s}_{\infty,\infty}}.\notag
\end{align} This means that the condition \eqref{+3} implies \eqref{+1}.
 This completes the proof.
\qquad $\hfill\Box$

\emph{\bf{Proof of Theorem \ref{main3}.}} If $s=\frac{1}{2}$, we see that $\theta=\frac{4s-2}{3-2s}=0$. The only difference between the case $s=\frac{1}{2}$ and the case $s>\frac{1}{2}$ comes from the term $J_1$ in \eqref {es24}. When $s=\frac{1}{2}$, we decompose $J_1$ as follows
\begin{align}\label{+11}
J_1=&\int_{\mathbb{R}^3}
\sum_{l=k-1}^{k+2}\sum_{l'=l-2}^{k-1}\dot{S}_{l-2}\dot{S}_{k}u_j\partial_j\dot{\Delta}_{l'}u_i\dot{\Delta}_lu^k_idx\\
=&\int_{\mathbb{R}^3}
\sum_{l=k-1}^{k+2}\sum_{l'=l-2}^{k-1}\dot{S}_{l-2}\dot{S}_{0}u_j\partial_j\dot{\Delta}_{l'}u_i\dot{\Delta}_lu^k_idx\notag\\
&+\int_{\mathbb{R}^3}
\sum_{l=k-1}^{k+2}\sum_{l'=l-2}^{k-1}\sum_{l''=0}^{k-1}\dot{S}_{l-2}\dot{\Delta}_{l''}u_j\partial_j\dot{\Delta}_{l'}u_i\dot{\Delta}_lu^k_idx.\notag
\end{align}From Bernstein inequality \eqref{Ber1}, \eqref{+10} and the fact $\dot{H}^\frac{1}{2}(\mathbb{R}^3)\hookrightarrow L^3(\mathbb{R}^3)$, we deduce that
\begin{align*}
&\lim_{k\to\infty}|\int_{\mathbb{R}^3}
\sum_{l=k-1}^{k+2}\sum_{l'=l-2}^{k-1}\dot{S}_{l-2}\dot{S}_{0}u_j\partial_j\dot{\Delta}_{l'}u_i\dot{\Delta}_lu^k_idx|\\
\leq & \lim_{k\to\infty}C\sum_{l=k-1}^{k+2}\sum_{l'=l-2}^{k-1}2^{l'}||\dot{S}_{0}u_j||_{L^\infty}||\dot{\Delta}_{l'}u_i||_{L^2}||\dot{\Delta}_lu_i||_{L^2}\notag\\
\leq & \lim_{k\to\infty}C||\dot{S}_{0}u_j||_{L^3}\sum_{l=k-1}^{k+2}\sum_{l'=l-2}^{k-1}2^{l'}||\dot{\Delta}_{l'}u_i||_{L^2}||\dot{\Delta}_lu_i||_{L^2}\notag\\
\leq& \lim_{k\to\infty}C||u_j||_{\dot{H}^\frac{1}{2}}\sum_{l=k-1}^{k+2}\sum_{l'=l-2}^{k-1}2^{\frac{l'}{2}}||\dot{\Delta}_{l'}u_i||_{L^2}2^{\frac{l}{2}}||\dot{\Delta}_lu_i||_{L^2}.\notag
\end{align*} From the fact $||u||_{\dot{H}^\frac{1}{2}}<\infty$, it is easy to see that $$\lim_{k\to\infty}\sum_{l=k-1}^{k+2}\sum_{l'=l-2}^{k-1}2^{\frac{l'}{2}}||\dot{\Delta}_{l'}u_i||_{L^2}2^{\frac{l}{2}}||\dot{\Delta}_lu_i||_{L^2}=0,$$ we thus get
\begin{align}\label{+12}
\lim_{k\to\infty}|\int_{\mathbb{R}^3}
\sum_{l=k-1}^{k+2}\sum_{l'=l-2}^{k-1}\dot{S}_{l-2}\dot{S}_{0}u_j\partial_j\dot{\Delta}_{l'}u_i\dot{\Delta}_lu^k_idx|=0\end{align} Substituting \eqref{+11}-\eqref{+12} and \eqref{es27}-\eqref{es30} into \eqref{es24},
we conclude that
\begin{align}\label{es37}
&\int_{\mathbb{R}^3}|(-\Delta)^\frac{1}{4}u|^2dx\\=&\liminf_{k\to\infty}\{\int_{\mathbb{R}^3}
\sum_{l=k-1}^{k+2}\sum_{l'=l-2}^{k-1}\sum_{l''=0}^{k-1}\dot{S}_{l-2}\dot{\Delta}_{l''}u_j\partial_j\dot{\Delta}_{l'}u_i\dot{\Delta}_lu^k_idx\notag\\
&+\int_{\mathbb{R}^3}\sum_{l=k-3}^{k}\sum_{l''=[\frac{k}{2}]}^{k-1}\dot{\Delta}_l\dot{S}_ku_j\partial_j\dot{\Delta}_{l''}u_i\tilde{\dot{\Delta}}_lu^k_idx
+\int_{\mathbb{R}^3}\sum_{l''=[\frac{k}{2}]}^{k-1}\sum_{k-1\leq l}\dot{\Delta}_lu^k_i\tilde{\dot{\Delta}}_{l}u^k_j\partial_j\dot{\Delta}_{l''}u_idx\}.\notag
\end{align}We thus establish the formula \eqref{snc3}.

To get the Liouville type theorem, we apply Bernstein inequalities \eqref{Ber1}-\eqref{Ber2} and \eqref{+10} to obtain that
\begin{align*}
&|\int_{\mathbb{R}^3}
\sum_{l=k-1}^{k+2}\sum_{l'=l-2}^{k-1}\sum_{l''=0}^{k-1}\dot{S}_{l-2}\dot{\Delta}_{l''}u_j\partial_j\dot{\Delta}_{l'}u_i\dot{\Delta}_lu^k_idx|\\
\leq&||\sum_{l''=0}^{k-1}\dot{\Delta}_{l''}u_j||_{L^\infty}\sum_{l=k-1}^{k+2}\sum_{l'=l-2}^{k-1}2^{l'}||\dot{\Delta}_{l'}u_i||_{L^2}||\dot{\Delta}_lu^k_i||_{L^2}\\
\leq& C||\sum_{l''=0}^{k-1}\dot{\Delta}_{l''}u_j||_{L^\infty}\sum_{l=k-1}^{k+2}\sum_{l'=l-2}^{k-1}2^{\frac{l'}{2}}||\dot{\Delta}_{l'}u_i||_{L^2}2^\frac{l}{2}||\dot{\Delta}_lu_i||_{L^2},
\end{align*}
\begin{align*}
&|\int_{\mathbb{R}^3}\sum_{l=k-3}^{k}\sum_{l''=[\frac{k}{2}]}^{k-1}\dot{\Delta}_l\dot{S}_ku_j\partial_j\dot{\Delta}_{l''}u_i\tilde{\dot{\Delta}}_lu^k_idx|\\
\leq &||\partial_j\sum_{l''=[\frac{k}{2}]}^{k-1}\dot{\Delta}_{l''}u_i||_{L^\infty}\sum_{l=k-3}^{k}||\dot{\Delta}_lu_j||_{L^2}||\tilde{\dot{\Delta}}_lu_i||_{L^2}\\
\leq& C2^k||\sum_{l''=[\frac{k}{2}]}^{k-1}\dot{\Delta}_{l''}u_i||_{L^\infty}\sum_{l=k-3}^{k}2^{-l}2^\frac{l}{2}||\dot{\Delta}_lu_j||_{L^2}2^{\frac{l}{2}}||\tilde{\dot{\Delta}}_lu_i||_{L^2}\\
\leq& C||\sum_{l''=[\frac{k}{2}]}^{k-1}\dot{\Delta}_{l''}u_i||_{L^\infty}\sum_{l=k-3}^{k}2^\frac{l}{2}||\dot{\Delta}_lu_j||_{L^2}2^{\frac{l}{2}}||\tilde{\dot{\Delta}}_lu_i||_{L^2},
\end{align*}
\begin{align*}
&|\int_{\mathbb{R}^3}\sum_{l''=[\frac{k}{2}]}^{k-1}\sum_{k-1\leq l}\dot{\Delta}_lu^k_i\tilde{\dot{\Delta}}_{l}u^k_j\partial_j\dot{\Delta}_{l''}u_idx|\\
\leq&||\partial_j\sum_{l''=[\frac{k}{2}]}^{k-1}\dot{\Delta}_{l''}u_i||_{L^\infty}\sum_{l\geq k-1}||\dot{\Delta}_lu_i||_{L^2}||\tilde{\dot{\Delta}}_{l}u_j||_{L^2}\\
\leq &C2^k||\sum_{l''=[\frac{k}{2}]}^{k-1}\dot{\Delta}_{l''}u_i||_{L^\infty}\sum_{l\geq k-1}2^{-l}2^\frac{l}{2}||\dot{\Delta}_lu_i||_{L^2}2^\frac{l}{2}||\tilde{\dot{\Delta}}_{l}u_j||_{L^2}\\
\leq& C||\sum_{l''=[\frac{k}{2}]}^{k-1}\dot{\Delta}_{l''}u_i||_{L^\infty}\sum_{l\geq k-1}2^\frac{l}{2}||\dot{\Delta}_lu_i||_{L^2}2^\frac{l}{2}||\tilde{\dot{\Delta}}_{l}u_j||_{L^2}.
\end{align*}Substituting the above estimates into \eqref{es37}, we have
\begin{align}\label{es38}
&\int_{\mathbb{R}^3}|(-\Delta)^\frac{1}{4}u|^2dx\\
\leq& C\liminf_{k\to\infty}(||\sum_{l''=0}^{k-1}\dot{\Delta}_{l''}u_j||_{L^\infty}+||\sum_{l''=[\frac{k}{2}]}^{k-1}\dot{\Delta}_{l''}u_j||_{L^\infty})\times\lim_{k\to\infty}(\sum_{l=k-1}^{k+2}\sum_{l'=l-2}^{k-1}2^{\frac{l'}{2}}
||\dot{\Delta}_{l'}u_i||_{L^2}2^\frac{l}{2}||\dot{\Delta}_lu_i||_{L^2}\notag\\
&+\sum_{l=k-3}^{k}2^\frac{l}{2}||\dot{\Delta}_lu_j||_{L^2}2^{\frac{l}{2}}||\tilde{\dot{\Delta}}_lu_i||_{L^2}+\sum_{l\geq k-1}2^\frac{l}{2}||\dot{\Delta}_lu_i||_{L^2}2^\frac{l}{2}||\tilde{\dot{\Delta}}_{l}u_j||_{L^2})\notag\\
\leq &C||u^0||_{L^\infty}\times\lim_{k\to\infty}(\sum_{l=k-1}^{k+2}\sum_{l'=l-2}^{k-1}2^{\frac{l'}{2}}
||\dot{\Delta}_{l'}u_i||_{L^2}2^\frac{l}{2}||\dot{\Delta}_lu_i||_{L^2}\notag\\
&+\sum_{l=k-3}^{k}2^\frac{l}{2}||\dot{\Delta}_lu_j||_{L^2}2^{\frac{l}{2}}||\tilde{\dot{\Delta}}_lu_i||_{L^2}+\sum_{l\geq k-1}2^\frac{l}{2}||\dot{\Delta}_lu_i||_{L^2}2^\frac{l}{2}||\tilde{\dot{\Delta}}_{l}u_j||_{L^2})\notag
\end{align}here $u^0= u-\dot{S}_0u$.
Since $||u||_{\dot{H}^\frac{1}{2}}<\infty$, it follows that
\begin{align}\label{es39}
&\lim_{k\to\infty}(\sum_{l=k-1}^{k+2}\sum_{l'=l-2}^{k-1}2^{\frac{l'}{2}}
||\dot{\Delta}_{l'}u_i||_{L^2}2^\frac{l}{2}||\dot{\Delta}_lu_i||_{L^2}
+\sum_{l=k-3}^{k}2^\frac{l}{2}||\dot{\Delta}_lu_j||_{L^2}2^{\frac{l}{2}}||\tilde{\dot{\Delta}}_lu_i||_{L^2}\\&+\sum_{l\geq k-1}2^\frac{l}{2}||\dot{\Delta}_lu_i||_{L^2}2^\frac{l}{2}||\tilde{\dot{\Delta}}_{l}u_j||_{L^2})=0.\notag
\end{align} Consequently, we deduce from \eqref{es38} and \eqref{es39} that the condition \eqref{+2} implies $$\int_{\mathbb{R}^3}|(-\Delta)^\frac{1}{4}u|^2dx=0$$ then $u\equiv0$. This completes the proof.
\qquad $\hfill\Box$

\section*{Acknowledgments}
The author is supported by the Construct Program of the Key Discipline in Hunan Province and NSFC (Grant No. 11871209), and the Hunan Provincial NSF (No. 2022JJ10033)
\\
\noindent {\bf Conflict of interest:}
We declare that we do not have any commercial or associative interest
that represents a conflict of interest in connection with the work
submitted.

\end{document}